\documentclass[a4paper,12pt]{jcpman}  

\usepackage[english]{babel}   
\usepackage[T1]{fontenc}     

\usepackage{mathrsfs}   
\usepackage{bm}        
\usepackage{color, xcolor}   

\usepackage{amsmath,amsfonts,amssymb}
\usepackage{extarrows}

\usepackage{algorithm}
\usepackage{algpseudocode}

\usepackage{siunitx}

\usepackage{graphicx}   
\usepackage{float}
\usepackage{tikz}       
\usepackage{palatino}
\usetikzlibrary{arrows,shapes,chains}

\usepackage{setspace}         
\usepackage{tocloft}          

\usepackage{subfig}     
\usepackage{booktabs}  
\usepackage{longtable} 
\graphicspath{{figures/}}

\usepackage{CJKutf8}

\usepackage{indentfirst}

\usepackage[justification=centering]{caption}
\makeatletter
\newcommand\dif{
  \mathop{}\!%
  \ifustc@math@style@TeX
    d%
  \else
    \mathrm{d}%
  \fi
}
\makeatother

\makeatletter
\renewcommand{\maketag@@@}[1]{\hbox{\m@th\normalsize\normalfont#1}}%
\makeatother

\usepackage{hyperref}  
\hypersetup{
    breaklinks = true,
    colorlinks = true,
    citecolor = {blue},
    urlcolor = {blue},
    linkcolor = {blue}
}

\usepackage{cleveref}

 \global\long\def\calK{\mathcal{K}}

 \global\long\def\calMh{\hat{\mathcal{M}}}

 \global\long\def\calRh{\hat{\mathcal{R}}}

 \global\long\def\scrG{\mathscr{G}}

 \global\long\def\scrO{\mathscr{O}}

 \global\long\def\scrSh{\hat{\mathscr{S}}}

 \global\long\def\scrvh{\rmvh_{ab}}

 \global\long\def\frakC{\mathfrak{C}}

 \global\long\def\frakCh{\hat{\mathfrak{C}}}

 \global\long\def\rmd{\mathrm{d}}
 \global\long\def\rme{\mathrm{e}}

 \global\long\def\rmv{v}

 \global\long\def\rmvh{\hat{{v}}}

 \global\long\def\bfe{\mathbf{e}}

 \global\long\def\bfu{\mathbf{u}}
 \global\long\def\bfv{\mathbf{v}}

 \global\long\def\bfuh{\hat{\mathbf{u}}}
 \global\long\def\bfvh{\hat{\mathbf{v}}}

 \global\long\def\bbN{\mathbb{N}}

 \global\long\def\bbR{\mathbb{R}}

 \global\long\def\nh{\hat{n}}
 \global\long\def\uh{\hat{u}}
 \global\long\def\Ih{\hat{I}}
 \global\long\def\Kh{\hat{K}}
 \global\long\def\Th{\hat{T}}
 \global\long\def\Mh{\hat{M}}
 \global\long\def\xih{\hat{\xi}}
 
\global\long\def\rmDelta{\mathrm{\Delta}}

 \global\long\def\vb{\bfv_b}
 \global\long\def\v{\bfv}

 \global\long\def\vhb{\bfvh_b}
 \global\long\def\vh{\bfvh}

 \global\long\def\vnvbth{\vh_{ab}}
 \global\long\def\vvbth{\rmvh_{ab}} 

 \global\long\def\rmvhb{\rmvh_b}
 
 \global\long\def\rmvA{\rmv_{\alpha}}
 \global\long\def\rmvhA{\rmvh_{\alpha}}

 \global\long\def\Pl{P_l}
 \global\long\def\PlI{P_l^1}
 \global\long\def\PlII{P_l^2}
 
 \global\long\def\PlIcotQ{P_l^{1,\mu}}
 \global\long\def\PlIlII{P_l^{1,2}}

 \global\long\def\lMI{l_{M_1}}

 \global\long\def\PL{P_L}
 \global\long\def\PLI{P_L^1}

 \global\long\def\PLIcotQ{P_L^{1,\mu}}
 \global\long\def\PLILII{P_L^{1,2}}

 \global\long\def\LMI{L_{M_1}}

 \global\long\def\mub{{\mu}_\beta}

 \global\long\def\wmub{w_{{\mu}_\beta}}

 \global\long\def\Dt{{\Delta}_t}
 \global\long\def\pDt{\delta_t}
 \global\long\def\pdt{\partial_t}

 \global\long\def\ddt{\frac{\partial}{\partial t}}

 \global\long\def\ddbfv{\nabla_{\v}}
 
 \global\long\def\ddbfvh{\nabla_{\vh}}

 \global\long\def\ddrmvh{\frac{\partial}{\partial \rmvh}}
 
 \global\long\def\dddrmvh{\frac{\partial^2}{\partial {\rmvh}^2}}
 
 \global\long\def\ddvnvbth{\nabla_{\vnvbth}} %
 \global\long\def\ddscrvh{\frac{\partial}{\partial \vvbth}}
 \global\long\def\dddscrvh{\frac{\partial^2}{\partial {\vvbth}^2}}

 \global\long\def\sumloq{\sum_{l=0}^{\infty} }
 \global\long\def\sumLoq{\sum_{L=0}^{\infty} }
 \global\long\def\sumlolM{\sum_{l=0}^{l_M} }
 \global\long\def\sumLoLM{\sum_{L=0}^{L_M} }
 \global\long\def\sumlolmax{\sum_{l=0}^{l_{max}} }
 \global\long\def\sumLoLmax{\sum_{L=0}^{l_{max}} }

 \global\long\def\sumLILMax{\sum_{L=1}^{l_{max}} }
 \global\long\def\sumlIlMax{\sum_{l=1}^{l_{max}} }

 \global\long\def\sumbIlMI{\sum_{\beta=1}^{\lMI} }
 \global\long\def\sumbILMI{\sum_{\beta=1}^{\LMI} }

 \global\long\def\e{\rme}

 \global\long\def\sinQ{\sqrt{1-\mu^2}}

 \global\long\def\Dvh{\rmDelta \rmvh}  
 \global\long\def\DvhA{\rmDelta \rmvhA}  

 \global\long\def\vthna{\frac{\rmv_{ath}^3}{n_a}}
 \global\long\def\navth{\frac{n_a}{\rmv_{ath}^3}}
 \global\long\def\nbvth{\frac{n_b}{\rmv_{bth}^3}}

 \global\long\def\fl{f_l}

\global\long\def\fh{\hat{f}}

 \global\long\def\fhl{\fh_l}
 \global\long\def\fho{\fh_0}

\global\long\def\fh{\hat{f}}

\global\long\def\Fh{\hat{F}}

 \global\long\def\fhl{\fh_l}
 \global\long\def\fho{\fh_0}

 \global\long\def\dfhl{\frac{\mathrm{\partial} \fhl}{\mathrm{\partial} \hat{\rmv}}}

\global\long\def\Fh{\hat{F}}
\global\long\def\Hh{\hat{H}}
\global\long\def\Gh{\hat{G}}

 \global\long\def\FhL{\hat{F}_L}
 \global\long\def\HhL{\hat{H}_L}
 \global\long\def\GhL{\hat{G}_L}
 \global\long\def\Fho{\Fh{_0}}
 \global\long\def\Hho{\Hh{_0}}
 \global\long\def\Gho{\Gh{_0}}

 \global\long\def\dGhL{\frac{\partial \GhL}{\partial \scrvh}}

\global\long\def\calMhjl{\calMh_{j,l}}

  \global\long\def\calMhII{{\calMh}{_{1,1}}}
  
  \global\long\def\calMhoo{{\calMh}{_{0,0}}}

\global\long\def\calRhjl{\calRh_{j,l}}

 \global\long\def\Kl{\calK_l}

 \global\long\def\KL{\calK_L}

 \global\long\def\IjFhL{I_{j,L}}

 \global\long\def\IiFhL{I_{i,L}}

 \global\long\def\ILFhL{I_{L,L}}
 \global\long\def\ILIIFhL{I_{L+2,L}}
 
 \global\long\def\JjFhL{J_{j,L}}

 \global\long\def\JiFhL{J_{i,L}}

 \global\long\def\JLpFhL{J_{L+1,L}}
 \global\long\def\JLnFhL{J_{L-1,L}}

 \global\long\def\lnAab{\ln{ \Lambda_{ab}}}

 \global\long\def\Gab{\Gamma_{ab}}

 \global\long\def\colhaa{\frakCh_{aa}}

 \global\long\def\cola{\frakC}
 
 \global\long\def\colla{{\frakC{_l}}}

 \global\long\def\colha{\frakCh}

 \global\long\def\colhla{{\frakCh{_l}}}

 \global\long\def\colab{\frakC_{ab}}

 \global\long\def\colhab{\frakCh_{ab}}
 
 \global\long\def\colhlab{{\frakCh{_l}}_{ab}}
 \global\long\def\colhoab{{\frakCh{_0}}_{ab}}
 \global\long\def\colhIab{{\frakCh{_1}}_{ab}}
 
 \global\long\def\colba{\frakC_{ba}}

 \global\long\def\colhoba{{\frakCh{_0}}_{ba}}
 \global\long\def\colhIba{{\frakCh{_1}}_{ba}}
 
 \global\long\def\Sh{\scrSh}

 \global\long\def\CHh{C_{\Hh}}
 \global\long\def\CGh{C_{\Gh}}

 \global\long\def\mp{m_p}

 \global\long\def\ua{\bfu_a}

\global\long\def\uha{\bfuh_a}

\global\long\def\uza{u_a}
\global\long\def\uzb{u_b}

\global\long\def\uzha{\uh_a}

\global\long\def\uzhas{\uh_{a_s}}

\global\long\def\uzhbr{\uh_{b_r}}

\global\long\def\Izha{\Ih_a}

 \global\long\def\vth{\rmv_{th}}
 \global\long\def\vath{\rmv_{ath}}
 \global\long\def\vabth{\rmv_{abth}}
 
 \global\long\def\veth{\rmv_{eth}}
 \global\long\def\vDth{\rmv_{Dth}}
 
 \global\long\def\vAth{\rmv_{\alpha th}}

 \global\long\def\vhath{\rmvh_{ath}}

 \global\long\def\vhaths{\rmvh_{{ath}_s}}

 \global\long\def\vbth{\rmv_{bth}}

 \global\long\def\vhbthr{\rmvh_{{bth}_r}}

 \global\long\def\ns{n_s}
 \global\long\def\Izs{I_s}
 
 \global\long\def\Ks{K_s}
 \global\long\def\ss{s_s}
 
\global\long\def\nha{\nh_a}
\global\long\def\nhb{\nh_b}

\global\long\def\nhas{\nh_{a_s}}

\global\long\def\nhao{\nh_{a_0}}
\global\long\def\nhaI{\nh_{a_1}}
\global\long\def\nhaII{\nh_{a_2}}

\global\long\def\nhbr{\nh_{b_r}}

\global\long\def\Kha{\Kh_a}
\global\long\def\Khb{\Kh_b}

\global\long\def\Tha{\Th_a}

\global\long\def\Rdtvath{\frac{1}{\rmv_{ath}} \dtvath}
\global\long\def\dtvath{\frac{\partial \rmv_{ath}}{\partial t}}

 \global\long\def\opo{\omega_{p_0}}  
 
 \global\long\def\nuTab{\nu_T^{ab}}

 \global\long\def\tauTab{\tau_T^{ab}}

 \global\long\def\NK{N_K}
 \global\long\def\NKo{N_{K_0}}
 \global\long\def\NKmax{N_K^{max}}
 \global\long\def\NKmin{N_K^{min}}
\global\long\def\NKa{N_{K_a}}
\global\long\def\NKb{N_{K_b}}

 \global\long\def\tko{t_0}

 \global\long\def\DtExp{{\Delta}_{\tko}^{Exp}}
 
 \global\long\def\tk{t_k}
 \global\long\def\Dtk{{\Delta}_{\tk}}

 \global\long\def\tkI{t_{k+1}}
 \global\long\def\DtkI{{\Delta}_{\tkI}}

 \global\long\def\Izsko{I_s^0}
 
 \global\long\def\Ksko{K_s^0}
 \global\long\def\ssko{s_s^0}

 \global\long\def\nako{n_a^0}
 
 \global\long\def\uzako{u_a^0}
 
 \global\long\def\Tako{T_a^0}

 \global\long\def\Izsk{I_s^k}
 
 \global\long\def\Ksk{K_s^k}
 \global\long\def\ssk{s_s^k}

 \global\long\def\nak{n_a^k}
 
 \global\long\def\uzak{u_a^k}
 \global\long\def\Izak{I_a^k}
 \global\long\def\Tak{T_a^k}
 \global\long\def\Kak{K_a^k}

 \global\long\def\vathk{\vath^k}

 \global\long\def\nkI{n^{k+1}}
 
 \global\long\def\IzkI{I^{k+1}}
 
 \global\long\def\KkI{K^{k+1}}

 \global\long\def\nakI{n_a^{k+1}}
 \global\long\def\rhoakI{\rho_a^{k+1}}
 \global\long\def\uzakI{u_a^{k+1}}
 \global\long\def\IzakI{I_a^{k+1}}
 
 \global\long\def\KakI{K_a^{k+1}}

 \global\long\def\vathkI{\vath^{k+1}}

 \global\long\def\uzhakI{\hat{u}_a^{k+1}}

 \global\long\def\navthkI{\frac{n_a^{k+1}}{\left(\vath^{k+1} \right)^3}}

 \global\long\def\Tbk{T_b^k}

 \global\long\def\IzbkI{I_b^{k+1}}
 
 \global\long\def\KbkI{K_b^{k+1}}

 \global\long\def\nbvthkI{\frac{n_b^{k+1}}{\left(\vbth^{k+1} \right)^3}}

 \global\long\def\vnvbth{\vh_{ab}}
 \global\long\def\vvbth{\rmvh_{ab}}

 \global\long\def\ddvnvbth{\nabla_{\vnvbth}} %

 \global\long\def\ddscrvh{\frac{\partial}{\partial \vvbth}}
 \global\long\def\dddscrvh{\frac{\partial^2}{\partial {\scrvh}^2}}
 
 \global\long\def\FIG#1{~\ref{#1}}
 \global\long\def\EQ#1{~(\ref{#1})}
 \global\long\def\EQo#1{(\ref{#1})}
 
 \global\long\def\EQgive#1{~(give in Eq.~(\ref{#1}))}
 \global\long\def\SEC#1{~\ref{#1}}
 \global\long\def\APP#1{~\ref{#1}}

\usepackage[title]{appendix}

\begin{document}

\title{A conservative, implicit solver for 0D-2V multi-species nonlinear Fokker-Planck collision equations}

\author[a]{Yanpeng Wang}
 
\author[a]{Jianyuan Xiao \footnote{xiaojy@ustc.edu.cn}}
\author[c]{Yifeng Zheng}
\author[d]{Zhihui Zou}
\author[b]{Pengfei Zhang}

\author[a]{Ge Zhuang}

\affil[a]{School of Nuclear Sciences and Technology, University of Science and Technology of China, Hefei, 230026, China}
\affil[b]{School of Physical Sciences, University of Science and Technology of China, Hefei, 230026, China}
\affil[c]{Institute of Plasma Physics, Chinese Academy of Sciences, Hefei 230031, China}
\affil[d]{Institute for Fusion Theory and Simulation, School of Physics, Zhejiang University, Hangzhou 310027, China}

\renewcommand{\cftdotsep}{\cftnodots}
\cftpagenumbersoff{figure}
\cftpagenumbersoff{table} 


\maketitle

\setlength{\parskip}{0pt}

\begin{abstract}
In this study, we present an optimal implicit algorithm specifically designed to accurately solve the multi-species nonlinear 0D-2V axisymmetric Fokker-Planck-Rosenbluth (FPR) collision equation while preserving mass, momentum, and energy. Our approach relies on the utilization of nonlinear Shkarofsky's formula of FPR (FPRS) collision operator in the spherical-polar coordinate. The key innovation lies in the introduction of a new function named King, with the adoption of the Legendre polynomial expansion for the angular coordinate and King function expansion for the speed coordinate. The Legendre polynomial expansion will converge exponentially and the King method, a moment convergence algorithm, could ensure the conservation with high precision in discrete form. Additionally, post-step projection onto manifolds is employed to exactly enforce symmetries of the collision operators. Through solving several typical problems across various nonequilibrium configurations, we demonstrate the high accuracy and superior performance of the presented algorithm for weakly anisotropic plasmas. 
\end{abstract}


\begin{spacing}{0.35}  

\section{Introduction}
\label{sect:Introduction} 

In plasma physics, the Fokker-Planck collision operator, known as the Fokker-Planck-Rosenbluth \cite{Rosenbluth1957,Taitano2015AEquation,Shkarofsky1963,Shkarofsky1967} (FPR)  or equivalently the Fokker–Planck–Landau \cite{Landau1937} (FPL) operator, is a fundamental tool for describing Coulomb collisions between particles under the assumptions of binary, grazing-angle collisions. This operator is particularly valuable for modeling various plasma systems, including those found in laboratory settings such as magnetic confinement fusion (MCF) and inertial confinement fusion (ICF), as well as in natural environments like Earth's magnetosphere and astrophysical phenomena like solar coronal plasma. When coupled with Vlasov's equation \cite{Vlasov1968THEGAS} and Maxwell's equations, it provides a comprehensive description of weakly coupled plasma across all collisionality regimes.

The FPR collision operator ensures strict conservation of mass, momentum and energy, while also adhering to the well-established H-theorem\cite{Boltzmann1872} , which guarantees that the entropy of the plasma system increases monotonically with time unless the system reaches a thermal equilibrium state. Throughout history, various formulations of the Fokker-Planck collision operator have been developed to suit different computational and theoretical needs. The FPL collision operator employs a direct integral formulation, making it ideal for conservative algorithms and the H-theorem \cite{Boltzmann1872} due to its symmetric nature. Conversely, the standard FPR collision operator \cite{Rosenbluth1957}  represents integral relationships using Rosenbluth potentials, which satisfies the Poisson equation in velocity space. The divergence form of the FPR (FPRD) \cite{Chang1970,Taitano2015AEquation,Taitano2016AnRegimes} collision operator is widely favored in numerical simulations due to its efficiency in fast solvers. Additionally, when employing spherical harmonic expansions, Shkarofsky’s formula of the FPR collision operator (FPRS) 
\cite{Shkarofsky1963,Shkarofsky1967} collision operator is often preferred for its computational advantages. These various formulations offer flexibility and efficiency in solving Vlasov-Fokker-Planck (VFP) equation, catering to different computational and theoretical requirements.
  
Historically, numerous efforts have been dedicated to addressing the numerical solution of the Fokker-Planck collision equation. Thomas\cite{Thomas2012} and Bell\cite{Bell2006FastEquation} reviewed the different numerical models of Fokker-Planck collision operator for ICF plasma. Cartesian tensor expansions\cite{Johnston1960,Shkarofsky1997,Kingham2004AnFields,Thomas2009} (CTE) and spherical harmonic expansions\cite{Bell2006FastEquation,Robinson2008ArtificialPulses,Tzoufras2011APhysics,Tzoufras2013APlasmas,Wu2013KineticIgnition,Mijin2020KineticLayer} (SHE), or Legendre polynomial expansions\cite{Krook1977,Bell1981,Matte1982,Shkarofsky1992,Alouani-Bibi2005,Zhao2018SimulationsWaves} when axisymmetric, are employed to handle the Fokker-Planck collision operator, which are considered equivalent to each other\cite{Johnston1960} .

SHE\cite{Bell2006FastEquation,Bell1981,Matte1982,Shkarofsky1992,Alouani-Bibi2005,Tzoufras2011APhysics} is an crucial method for moderate nonequilibrium plasma when the ratio of average velocity to thermal velocity is not large. As emphasized by Bell\cite{Bell2006FastEquation} , the amplitude of each harmonics will decay exponentially at a rate proportional to ${l(l+1)/2}$. Even in cases of weak collisions, this $l(l+1)/2$ leads to strong damping of higher-order spherical harmonics, naturally terminating the expansion. Early studies by Bell \cite{Bell1981} and Matte \cite{Matte1982} focused on including the first two order harmonics to investigate non-Spitzer heat flow in ICF plasma. Subsequent work by Shkarofsky \cite{Shkarofsky1992} and Alouani-Bibi \cite{Alouani-Bibi2005} extended this approach to higher orders, resulting in the widely used semi-anisotropic collision operator \cite{Alouani-Bibi2005} . In recent years, several VFP codes \cite{Bell2006FastEquation,Tzoufras2011APhysics,Wu2011RelativisticScenario} based on the semi-anisotropic model have been developed. However, effectively calculating the full nonlinear collision operator in the SHE approach remains a challenge \cite{Bell2006FastEquation,Tzoufras2011APhysics} , especially in scenarios involving large mass disparities such as electron-deuterium collisions in fusion plasma. While previous simulations utilizing SHE have adopted the semi-anisotropic model and maintained mass and energy conservation, achieving exact momentum conservation in discrete simulations remains problematic.

Other computational approaches such as meshfree methods\cite{Pareschi2000FastOperator,Filbet2002ACase,Pataki2011FastOperator,Askari2015MeshlessEquation} and finite volume method\cite{MortonK.W.2005NumericalEquations} (FVM) are also employed to solve the Fokker-Planck collision equation. Fast spectral method based on FFT\cite{Pareschi2000FastOperator,Filbet2002ACase,Pataki2011FastOperator} or Hermite polynomial expansion\cite{Li2021HermitePlasma} has shown rapid convergence of spectral expansion strategy\cite{Press2007NumericalRecipes} in approximating the FPL collision operator. Additionally, Askari\cite{Askari2015MeshlessEquation} employed a meshfree method using multi-quadric radial basis functions (RBFs) to approximate the solution of the 0D-1V Fokker-Planck collision equation. 
Taitano et al. \cite{Taitano2015AEquation,Taitano2015II,Taitano2016AnRegimes,Taitano2017AnGeometry,Daniel2020AElectrons} carried out systematic studies based on the 0D-2V FPRD collision operator by directly discretizing the collision equation with FVM. They\cite{Taitano2015AEquation}
developed an implicit algorithm and overcame the Courant-Friedrichs-Lewy\cite{Courant1986UberPhysik} (CFL) condition by utilizing a second-order BDF2 implicit solver and employing the multigrid (MG) method\cite{Saad2003IterativeEdition} in Jacobian-Free Newton-Krylov (JFNK)\cite{Saad2003IterativeEdition} solver. 
Furthermore, by normalizing the velocity space to the local thermal velocities of each species individually \cite{Larroche2003KineticImplosions} , works in Ref. \cite{Taitano2016AnRegimes} developed a discrete conservation strategy that utilizes nonlinear constraints to enforce the continuum symmetries of the collision operator. However, those approaches did not take advantages of Coulomb collisions, similar to SHE\cite{Bell2006FastEquation} , to reduce the number of meshgrids when there are no distinguishing asymmetries in the velocity space. 

The challenge of employing SHE \cite{Bell2006FastEquation} and previous meshfree\cite{Askari2015MeshlessEquation} approaches lies in embedding discrete conservation laws within the numerical scheme. According to manifold theory \cite{HairerE.2006GeometricIntegration} , maintaining a small local error through post-step projection onto manifolds preserves the same convergence rate. Therefore, backward error analysis\cite{Moore2003BackwardMethods,Reich1999BackwardIntegrators} has become a crucial tool for understanding the long-term behavior of numerical integration methods and preserving conservation properties in the numerical scheme.

In this study, our objective is to address the full nonlinearity, discrete conservation laws, and the temporal stiffness challenge of the 0D-2V axisymmetric multi-species FPRS collision equations within the SHE approach. Similar to previous works in Ref. \cite{Taitano2016AnRegimes} , we normalize the velocity spaces to the local thermal velocities for all species separately. However, instead of utilizing multigrid (MG) technology as in Ref. \cite{Taitano2015AEquation} , we employ a meshfree\cite{Liu2005} approach based on King method (details in Sec.\SEC{King method}) by introducing a novel shape function named King to overcome the classical CFL condition. To tackle the nonlinear, stiff FPRS collision equations, we propose an implicit algorithm based on Legendre polynomial expansion for the angular coordinate, the King function expansion for the speed coordinate, and the trapezoidal method  \cite{Press2007NumericalRecipes,Rackauckas2017DifferentialEquations.jlaJulia} for time integration. Romberg integration \cite{Bauer1961AlgorithmIntegration} is employed to compute kinetic moments with high precision, and backward error analysis \cite{Hairer1999BackwardMethods,Reich1999BackwardIntegrators,Moore2003BackwardMethods} is applied to ensure numerical conservation of mass, momentum, and energy. The H-theorem \cite{Boltzmann1872} is satisfied in the discretization scheme and utilized as a criterion for convergence of our algorithm.

The rest of this paper is organized as follows. Sec.\ref{The FPR collision equation} introduces the FPRS collision equation and its normalization. The discretization of the nonlinear FPRS collision equation is given in Sec.\ref{Discretization of the nonlinear FPRS collision equation}, encompassing angular discretization and the King method for the speed coordinate. An implicit time discretization and conservation strategies is discussed in Sec.\ref{Time integration}. The numerical performance of our solver, including accuracy and efficiency, is demonstrated with various multi-species tests in Sec.\ref{Numerical results}. Finally, we conclude our work in Sec.\ref{Conclusion}.

\end{spacing}

\begin{spacing}{0.7}   

\section{The Fokker-Planck-Rosenbluth collision equation}
\label{The FPR collision equation}

 The relaxation of Coulomb collision in a spatially homogeneous multi-species plasma can be described by the FPR collision equation. For species $a$, the velocity distribution functions $f(\v,t)$, in velocity space $\v$, satisfies:
  \begin{eqnarray}
      \ddt f \left(\v,t \right) &=& \cola \left(\v,t \right) ~.\label{FPeq}
  \end{eqnarray}
 The term on the right-hand side represents the FPR collision operator of species $a$. In this paper, the mass, time, charge, thermal velocity, number density $n_a$, temperature $T_a$ and permittivity are normalized by the proton mass $\mp$, characteristic time $\tau_0$, elementary charge $|q_e|$, vacuum speed of light $c_0$, reference number density $n_0=10^{20} \mathrm m^{-3}$, practical unit $T_k=\mathrm{keV} $ and permittivity of vacuum $\varepsilon_0$ respectively. The dimensionless form of other quantities is determined based on their correlation with these dimensionless quantities. 
 
 The normalization FPR collision equation maintains the same structure as Eq.\EQ{FPeq} and $\cola$ can be formulated as: 
  \begin{eqnarray}
      \cola \left(\v,t \right) &=& \sum_{b=1}^{N_s} \colab ~.\label{cola}
  \end{eqnarray}
  Here, $N_s$ represents the total number of plasma species and function $\colab$ denotes the 
  FPRS collision operator\cite{Shkarofsky1963,Shkarofsky1967} for species $a$ colliding with species $b$, given by:
  \begin{eqnarray}
      \colab \left(\v,t \right) &=& \Gab \left [4 \pi m_M F f + \left(1-m_M \right) \ddbfv H \cdot \ddbfv f + \frac{1}{2} \ddbfv \ddbfv G : \ddbfv \ddbfv f \right]  ~. \label{FPRS}
  \end{eqnarray}
  The mass ratio is denoted as $m_M=m_a/m_b$, where symbols $m_a$ and $m_b$ represent the masses of species $a$ and $b$ respectively. Parameter $\Gab=C_{\Gamma} \times 4\pi \left(\frac{Z_a Z_b}{4\pi m_a} \right)^2 \lnAab$ and the dimensionless coefficient $ C_{\Gamma}= {\tau_0 \opo^4} / {(n_0 c_0^3)} $ where $\opo = \sqrt{n_0 q_e^2 / (m_p \varepsilon_0)}$. Symbols $Z_a$ and $Z_b$ denote the ionization state of species $a$ and $b$. Parameter $\lnAab$ represents the Coulomb logarithm\cite{Huba2011NRLFORMULARY} of species $a$ and $b$, which is a weak function of the number of particles in the Debye sphere. Function $F=F(\vb,t)$, representing the distribution function of background species $b$. Functions $H$ and $G$ denote the Rosenbluth potentials, which are integral operators for the background distribution function $F$, reads:
  \begin{eqnarray}
    H(\v) &=& \int \frac{1}{\left|\v-\vb \right|} F\left(\vb,t \right) \rmd \v _b, \label{H}
    \quad 
    G(\v) \ = \ \int \left|\v-\vb \right| F\left(\vb,t \right) \rmd \v _b ~. \label{G}
  \end{eqnarray}
  
  It is worth noting that, in order to reduce the burden of subscripts in the subsequent sections, the quantities of species $a$, such as the distribution function $f$ and velocity $\v$, will not include the subscript $a$ in this paper, as commonly practiced\cite{wang2024Relaxationmodel,Bell2006FastEquation,Tzoufras2013APlasmas,Matte1982,Taitano2015AEquation,Taitano2015II,Shkarofsky1963,Shkarofsky1967,Shkarofsky1992} . Similarly, the quantities of species $b$, denoted by capital letters, such as the distribution function $F$, Rosenbluth potentials $H$ and $G$, are also omitted from including the subscript $b$.

\subsection{Conservation}
\label{Conservation}

  The FPRS collision operator\EQgive{FPRS} preserves mass, momentum and energy which stems from its symmetries\cite{Braginskii1965TransportPlasma} . Using the inner product definition $ \left<f,g \right>_{\v}= \int g(\v) f(\v) \rmd \v$, these conservation laws can be expressed as follows:
  \begin{eqnarray}
      \left<1, \colab \right>_{\v} &=& \left<1, \colba \right>_{\v} \quad \equiv \quad 0, \label{Cn}
      \\
      m_a \left<\v, \colab \right>_{\v} &=& - m_b \left<\v, \colba \right>_{\v}, \label{CI}
      \\
      m_a \left<\frac{\v^2}{2}, \colab \right>_{\v} &=& - m_b \left<\frac{\v^2}{2}, \colba \right>_{\v} ~. \label{CK}
  \end{eqnarray}
  
  In theory, the FPRS collision operator satisfies the well-known H-theorem. By defining the Boltzmann's entropy of species $a$, $s_a (t) = - \left<f, \ln{f}\right>_{\v}$, the total entropy of the plasma system can be expressed as $s_s(t)= \sum_a s_a$. According to the H-theorem, the total entropy of an isolated plasma system will monotonically increase over time unless there is no change in total entropy, indicating that all distribution functions are Maxwellian with a common temperature and average velocity.

  \begin{figure}[htp]
	\begin{center}
		\includegraphics[width=0.6\linewidth]{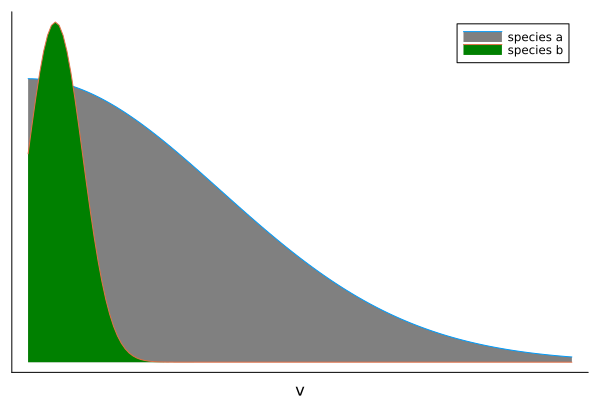}
	\end{center}
	\caption{Illustration of the velocity distribution functions in speed coordinate for disparate thermal velocities in a subsonic plasma system.}
	\label{FigfaFb}
  \end{figure}
  
\subsection{Normalization}
\label{Normalization}

In fusion plasma, the presence of disparate thermal velocities poses additional challenges, arising from the significant
mass discrepancy ($electron$-s$ion$ collisions) or energy difference ($deuterium$-$alpha$ collisions). This paper specifically focuses on the weakly anisotropic plasmas where the system exhibits no distinguishing asymmetries in the velocity space, such as when the average velocity is significantly smaller than the thermal velocity. Fig.\FIG{FigfaFb} illustrates the distribution function in speed coordinate when bulk flow of plasma is subsonic with a significant discrepancy in thermal velocities. This discrepancy adds complexity to discretizing the speed, particularly when mapping the background species distribution function to the collision species domain. Previous studies\cite{Larroche2003KineticImplosions,Taitano2016AnRegimes} have shown that normalizing the distribution function by its local thermal velocity (denoted as $\vath$ for species $a$) can help to mitigate these challenges. This normalization not only alleviates the need for different meshing requirements between multiple species but also ensures consistent evolution of thermal velocities over time with respect to temperature changes in distribution functions.
  
Therefore, we normalize the velocity space with its corresponding local thermal velocity, for species $a$, reads:
  \begin{equation}
      \vh = \v / \vath ~. \label{vh}
  \end{equation}
Thus, we have $\ddbfvh = \vath \ddbfv$ and the distribution function can be normalized as follow:
  \begin{eqnarray}
      \fh \left(\vh,t \right) &=& n_a^{-1} \vath^3 f \left(\v / \vath,t \right) ~. \label{fh}
  \end{eqnarray}
Therefore, the normalized background distribution function and Rosenbluth potentials can be written as:
  \begin{eqnarray}
      \Fh \left(\vhb,t \right) &=& n_b^{-1} \vbth^3 F \left(\vb / \vbth,t \right),  \label{Fh}
      \\
      \Hh(\vnvbth, t) &=& \frac{\vbth}{n_b }  H(\v / \vbth, t) \ = \ \int \frac{1}{\vnvbth -\vhb} \Fh \left(\vhb,t \right) \rmd \vhb, \label{Hh} 
      \\
      \Gh(\vnvbth, t) &=& \frac{1}{n_b \vbth} G(\v / \vbth, t) \ = \ \int \left(\vnvbth -\vhb \right) \Fh \left(\vhb,t \right) \rmd \vhb, \label{Gh}
  \end{eqnarray}
where $\vhb=\vb / \vbth $ and $\vnvbth=\v / \vbth $.

Hence, the FPRS collision operator\EQgive{cola} of species $a$ can be normalized as follows:
  \begin{eqnarray}
      \colha \left(\vh,t \right) &=& \vthna \cola(\v / \vath, t) = \sum{_{b=1}^{N_s}} \nbvth \Gab \colhab, \label{FPShda}
  \end{eqnarray}
where
  \begin{eqnarray}
      \colhab \left(\vh,t \right) &=& 4\pi m_M \Fh(\vnvbth,t) \fh + \CHh \ddvnvbth \Hh \cdot \ddbfvh \fh + \CGh \ddvnvbth \ddvnvbth \Gh : \ddbfvh \ddbfvh \fh ~. \label{FPShdab}
  \end{eqnarray}
Here, $m_M$ is the mass ratio, $\ddbfvh$ and $\ddvnvbth$ are gradients in normalized velocity space $\vh$ and $\vnvbth$ respectively.
The coefficients in Eq.\EQ{FPShdab} are given by:
  \begin{eqnarray} 
      \CHh=\left(1-m_M \right) \vbth / \vath, \label{CHh}
      \quad 
      \CGh=\left(\vbth / \vath \right)^2 / 2
      ~. \label{CGh}        
  \end{eqnarray}
The concrete formulation of Eq.\EQ{FPShdab} is provided in Appendix\APP{Normalized FPRS collision operator}. The normalized like-particle collision operator can be derived from Eq.\EQ{FPShdab} by replacing $\Fh$ and $\vnvbth$ by $\fh$ and $\vh$ respectively, reads:
  \begin{eqnarray}
      \colhaa \left(\vh,t \right) &=& 4 \pi \fh \fh + \frac{1}{2} \ddbfvh \ddbfvh \Gh : \ddbfvh \ddbfvh \fh ~. \label{colhaa}
  \end{eqnarray}
  
After the applications of Eqs.\EQ{vh}-\EQo{fh} and Eq.\EQ{FPShda}, the finial form of FPR collision equation to be solved numerically in this study is:
  \begin{eqnarray}
      \ddt f \left(\vh,t \right) &=& \cola \left(\vh,t \right) , \label{FPeqd}
  \end{eqnarray}
where
  \begin{eqnarray}
      f \left(\vh,t \right) &=& \navth \fh \left(\vh,t \right),\label{f}
      \\
       \cola \left(\vh,t \right) &=& \navth \sum{_{b=1}^{N_s}} \nbvth \Gab \colhab ~. \label{Ca}
  \end{eqnarray}
This paper focuses on the scenario of axisymmetric velocity space. Therefore, Eq.\EQ{FPeqd} can be referred to as the 0D-2V FPRS collision equation.  It is important to note that the fully normalized form of both the velocity space and the objective function, such as $\fh (\vh)$, is utilized to derive $\cola \left(\vh,t \right)$, while the time derivative employs a semi-normalized form only on the velocity space, such as $f \left(\vh,t \right)$, rather than $\fh (\vh, t)$ used by Taitano\cite{Taitano2016AnRegimes} .

The utilization of this semi-normalized equation\EQ{FPeqd} can eliminate the non-inertial terms\cite{Taitano2016AnRegimes} and simplify the complexity of the FPRS collision equation. The successful implementation of this approach relies on the utilization of the King method introduced in Sec.\SEC{King method},  in combination with the time block technique (TBT) detailed in Sec. \ref{Time block technique}. The specific implementation in implicit algorithms is provided in the main procedure described by Algorithm \ref{alg: Trapez}. In order to develop an effective algorithm, we make the assumption that the distribution function, $f(\vh,t)$, is a smooth function in the velocity space. It is reasonable that the Coulomb collisions always tends to eliminate the fine structures of the distribution function\cite{Bell2006FastEquation} .

\section{Discretization of the nonlinear FPRS collision equation}
\label{Discretization of the nonlinear FPRS collision equation}

In this study, we adopt a meshfree approach\cite{Liu2005} (requiring field nodes) for discretizing the 0D-2V axisymmetric FPRS collision operator within a spherical-polar coordinate in velocity space. The SHE\cite{Johnston1960,Robson1986VelocityEquation} method is employed to discretize the angular coordinate of the distribution function. Subsequently, the King function expansion (KFE) method, presented in Sec.\SEC{King method}, is utilized to the speed coordinate. Based on this framework of SHE together with KFE, a moment convergence algorithm is developed for solving the 0D-2V FPRS collision equation represented by Eq.\EQ{FPeqd}.

\subsection{Angular discretization}
\label{Angular discretization}

We opt for the SHE\cite{Johnston1960,Robson1986VelocityEquation} method (utilizing Legendre polynomial expansions\cite{Rosenbluth1957,Matte1982,Shkarofsky1992,Sunahara2003Time-DependentImplosions,Alouani-Bibi2005,Zhao2018SimulationsWaves} when the velocity space is axisymmetric) to analytically adapt the velocity-space discretization in angular coordinate. Unlike previous studies\cite{Bell2006FastEquation, Tzoufras2013APlasmas, Joglekar2018ValidationPhysics} employing the semi-anisotropic model, we will maintain the full nonlinearity of the FPRS collision operator for all species.

\subsubsection{Legendre polynomial expansions}
\label{Legendre polynomial expansions}

The normalized distribution function of species $a$ is described by the real function $\fh(\rmvh, \theta,t)$ in axisymmetric systems, which can be expanded using Legendre polynomials in normalized velocity space $\vh (\rmvh, \theta)$. It can be expressed as:
  \begin{eqnarray}
      \fh \left(\rmvh,\mu,t \right) &=& \sumloq \fhl \left( \rmvh,t \right) \Pl \left(\mu \right) ~. \label{fhfhl}
  \end{eqnarray}
Here, $\rmvh=|\v| / \vath $, $\mu=\cos{\theta}$ and $0 \le \theta < \pi$ when choosing the symmetric axis to be the direction $z$ with base vector $\bfe_z$. The function $\Pl (\mu)$ represents the $l^{th}$-order Legendre polynomials. For the sake of convenience, we will henceforth denote Eq.\EQ{fhfhl} as the SHE in axisymmetric velocity space. The calculation for the $l^{th}$-order normalized amplitude $\fhl \left( \rmvh,t \right)$ can be obtained through the inverse transformation of Eq.\EQ{fhfhl} as follows:
   \begin{eqnarray}
       \fhl \left(\rmvh, t \right) &=& \frac{1}{(N_l)^2} \int_{-1}^1 \fh \left (\rmvh, \mu, t \right)   \Pl (\mu) \rmd \mu, \label{fhlfh}
   \end{eqnarray}
where $N_l=\sqrt{(2l+1)/(4 \pi)}$, representing the normalization coefficient in spherical harmonics\cite{Arfken1971MathematicalEdition} .

  \begin{figure}[H]
	\begin{center}
		\includegraphics[width=0.67\linewidth]{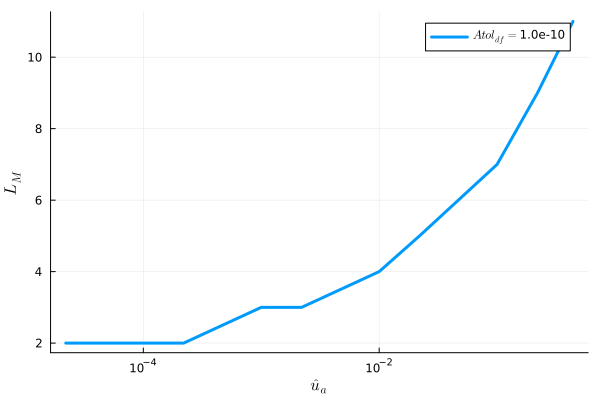}
	\end{center}
	\caption{Convergence of SHE for drift-Maxwellian distribution: $l_M$ as a function of $\hat{u}_a$ when $Atol_{df}=10^{-10}$.}
	\label{FigLMuh}
  \end{figure}

Due to exponential decay of the $l^{th}$-order amplitude at a rate proportion to $l(l+1)/2$, as stated in Ref.\cite{Bell2006FastEquation} , there is a natural termination to the expansions. Thus, the function $\fh(\rmvh, \mu,t)$ is represented by a finite set of amplitudes $\fhl \left( \rmvh,t \right)$, which are dependent on time $t$ and normalized speed $\rmvh \equiv |\vh|$. The series in Eq.\EQ{fhfhl} can be truncated at a maximum order, denoted as $l_M$, under the specified condition that $max|\fhl(\rmvh,t)|\le Atol_{df}$. Here, $l_M \in \bbN$ and the subscript $M$ is short for “max”. The convergence of SHE for drift-Maxwellian distribution is depicted in Fig.\FIG{FigLMuh}. It shows that $l_M$ is monotonic function of $\hat{u}$ when $\hat{u} \le 0.5$. For example with $\uzha \le 2.2 \times 10^{-3}$ and $Atol_{df}=10^{-10}$, it leads to a maximum order, $l_M = 3$.

Eq.\EQ{fhlfh} can be computed by utilizing Gaussian quadrature\cite{Press2007NumericalRecipes} ,
which can be reformulated as:
   \begin{eqnarray}
       \fhl \left(\rmvh, t \right) &=& \sumbIlMI \wmub  \Pl (\mub) \fh \left (\rmvh, \mub, t \right) + \scrO(Atol_{df}), \quad l = 0,1,2,\cdots,l_M ~, \label{fhlfhGQ}
   \end{eqnarray}
where 
   \begin{eqnarray}
       \lMI &=& l_M+1 ~. \label{lM1}
   \end{eqnarray}
The Gauss-Legendre abscissas will used as the field nodes for the polar angle coordinate $\mu$. The node $\mub$ represents the $\beta^{th}$ roots of the Legendre polynomial $\Pl(\mu)$, out of a total of $\lMI$ roots. The associated weight $\wmub$ is computed using Fornberg's algorithm \cite{Fornberg1998CalculationFormulas} .

\subsubsection{Rosenbluth potentials}
\label{Rosenbluth potentials}

Similar to Eq.\EQ{fhfhl}, the normalized distribution function of the background species $b$, can be expended as follows:
  \begin{eqnarray}
      \Fh \left(\rmvh_b,\mu,t \right) &=& \sumLoq \FhL \left( \rmvh_b,t \right) \PL \left(\mu \right) ~. \label{FhFhL}
  \end{eqnarray}With a maximum truncated order, $L_M$, the $L^{th}$-order normalized amplitude can be formulated as:
   \begin{eqnarray}
       \FhL \left(\rmvhb, t \right) &=& \sumbILMI \wmub \PL (\mub) \Fh \left (\rmvhb, \mub, t \right) + \scrO (Atol_{df}) , \label{FhLFhGQ}
   \end{eqnarray}
where $\LMI = L_M + 1$.
Here, we utilize $L$ rather than $l$ because $\Fh$ and $\fh$ typically exhibit distinct convergence rates in angular coordinate. The maximum order of angular coordinate during collisions between species $a$ and species $b$ will be denoted as $l_{max} = \max{(l_M,L_M)}$. Additionally, we will disregard the truncation error terms from this point onward.

In a similar manner, the normalized Rosenbluth potentials of species $a$, as described by equations\EQ{Hh}-\EQo{Gh} due to the presence of background species $b$, can also be represented in an expanded form, reads:
  \begin{eqnarray}
      \Hh \left(\vvbth,\mu ,t \right) &=& 4 \pi \sumLoLM \HhL \left(\vvbth ,t \right) \PL \left(\mu \right) , \label{HhHhL}
      \\
      \Gh \left(\vvbth,\mu,t \right) &=& 4 \pi \sumLoLM\GhL \left(\vvbth ,t \right) \PL \left(\mu \right) ~. \label{GhGhL}
  \end{eqnarray}
Here, $\vvbth= |\v| / \vbth $. The coefficient $4 \pi$ stems from the definitions represented by Eqs.\EQ{IjFL}-\EQo{JjFL}. The $L^{th}$-order amplitudes of $\Hh$ and $\Gh$ can be computed in the following integral form:
  \begin{eqnarray}
      \HhL \left(\vvbth,t \right) &=& \frac{1}{2 L + 1} \frac{1}{\vvbth} \left (\ILFhL + \JLpFhL \right) , \label{HhL}
      \\
      \GhL \left(\vvbth,t \right) &=& \frac{1}{2 L + 1} \frac{1}{\vvbth} \left (\frac{\ILIIFhL + \JLpFhL}{2 L + 3} - \frac{\ILFhL + \JLnFhL}{2 L - 1} \right)  ~. \label{GhL}
  \end{eqnarray}
Here, $\IiFhL$ and $\JiFhL$ represent the functionals of the normalized background distribution function $\FhL\left(\rmvhb,t \right)$, following a similar approach as Shkarofsky et al.\cite{Shkarofsky1967,Shkarofsky1997} , reads:
  \begin{eqnarray}
      \IiFhL \left(\vvbth,t \right) &=& I_i \left[\FhL \right] \ = \ \frac{1}{(\vvbth)^i} \int_0^{\vvbth} (\rmvhb)^{i+2} \FhL \left( \rmvhb,t \right) \rmd \rmvhb, \quad i = L, L + 2, \label{IjFL}
      \\ 
      \JiFhL \left(\vvbth,t \right) &=& J_i \left[\FhL \right] \ = \ (\vvbth)^i \int_{\vvbth}^{\infty} \frac{\rmvhb^2}{(\rmvhb)^i} \FhL \left( \rmvhb,t \right) \rmd \rmvhb, \quad i = L - 1, L + 1 ~. \label{JjFL}
  \end{eqnarray}
The aforementioned definitions do not include the coefficient $4 \pi$, which arises from the use of spherical-polar coordinates in velocity space. The Jacobian $(\rmvh_b)^2$ used in speed integrals, such as Eqs.\EQ{IjFL}-\EQo{JjFL}, and the subsequent definition of kinetic moment represented by Eq.\EQ{Mhjl}, also stems from the application of spherical-polar coordinates in velocity space.

\subsubsection{FPRS collision spectrum equation}
\label{FPRS collision spectrum equation}

Similar to the expansion of $\fh(\vh$)\EQgive{fhfhl}, the normalized FPRS collision operator\EQ{FPShda} can also be expanded based on the Legendre polynomials, reads:
  \begin{eqnarray} 
      \colha \left(\vh, t \right) &=& \sumlolM \colhla \left( \rmvh ,t \right) \Pl \left(\mu \right) ~. \label{FPShdl}
  \end{eqnarray}
The $l^{th}$-order amplitude of normalized multi-species nonlinear FPRS operator of species $a$ is given by:
  \begin{eqnarray}
      \colhla \left(\rmvh,t \right) &=& \sum_{b=1}^{N_s} \nbvth \Gab \colhlab ~. \label{FPShlda}
  \end{eqnarray}
Function $\colhlab$ represents the $l^{th}$-order amplitude of normalized FPRS collision operator for species $a$ colliding with species $b$, which can be expressed as:
   \begin{eqnarray}
       \colhlab \left(\rmvh, t \right) &=& \sumbIlMI \wmub \Pl (\mub) \colhab \left (\rmvh, \mub, t \right)
       ~. \label{collab}
   \end{eqnarray}
 Function $\colhab$ depends on $\fhl$ and $\FhL$, as detailed in Appendix\APP{Normalized FPRS collision operator}.
   
Substituting Eq.\EQ{fhfhl} and Eq.\EQ{FPShdl} into the FPRS collision equation\EQ{FPeqd} yields the following spectrum equation:
  \begin{eqnarray}
      \ddt \fl \left(\rmvh,t \right) &=& \colla \left(\rmvh,t \right), \quad 0 \le l \le l_M, \label{dtfl}
  \end{eqnarray}
where
  \begin{eqnarray}
      \fl \left(\rmvh,t \right) &=& \navth \fhl \left(\rmvh,t \right), \label{fl}
      \\
      \colla \left(\rmvh,t \right) &=& \navth \colhla \left(\rmvh,t \right)~. \label{Cl}
  \end{eqnarray}
Referring to above equation as the 0D-2V FPRS collision spectral equation.

\subsubsection{Moment constraints}

The $(j,l)^{th}$-order normalized kinetic moment of species $a$, denoted as $\Mh_{j,l}$, is generally defined as:
\begin{equation}
    \Mh_{j,l} \left(t \right) \ = \ \Mh_j \left[\fhl \right] \ = \ 4 \pi \int_0^{\infty} (\rmvh)^{j+2} \fhl \left(\rmvh,t \right) \rmd \rmvh ~. \label{Mhjl}
\end{equation}
The first few orders of $\Mh_{j,l}$ specifically relative to the conserved moments satisfy the following relations:
  \begin{eqnarray}
      \nha \left(t \right) &=& \Mh_{0,0},   \label{nha} 
      \\
      \Izha \left(t \right) &=& \frac{I_a}{\rho_a \vath} \ = \ \frac{1}{3} \Mh_{1,1},  \label{Iha} 
      \\
      \Kha \left(t \right) &=& \frac{K_a}{n_a T_a} \ = \ \Mh_{2,0} ~. \label{Kha}
  \end{eqnarray}
Theoretically, $\nha(t)$ is conserved and equal to 1. The momentum is defined as $I_a(t) = \rho_a (\ua \cdot \bfe_z)$ where the average velocity $\ua(t) =  n_a^{-1} \int \v f \left(\v,t \right) \rmd \v$, temperature $T_a(t) = m_a \vath^2 / 2$ and total energy $K_a(t) = \frac{m_a}{2} \int \v^2 f \left(\v,t \right) \rmd \v$. The normalized average velocity $\uzha(t) = (\ua \cdot \bfe_z)/\vath$, which is equivalent to $\Izha$. Similar to the normalized kinetic moment\EQ{Mhjl}, we define the $(j,l)^{th}$-order normalized kinetic dissipative force of species $a$ as follows:
\begin{equation}
    \calRhjl \left(t \right) \ = \ \calRh_j \left[\colhla \right] \ = \ 4 \pi \int_0^{\infty} (\rmvh)^{j+2} \colhla \left(\rmvh,t \right) \rmd \rmvh ~. \label{Rhjl}
\end{equation}

Multiplying both sides of the FPRS collision spectral equation\EQ{dtfl} by $4 \pi m_a \rmv^{j+2} \rmd \rmv$ and integrating over the semi-infinite interval $\rmv = \left [0, \infty \right)$, simplifying the result gives the $(j,l)^{th}$-order FPRS collision spectral equation in weak form:
\begin{eqnarray}
    \ddt \left[ 4 \pi \rho_a (\vath)^j \int_0^\infty (\rmvh)^{j+2} \fhl (\rmvh,t) \rmd \rmvh \right] &=& 4 \pi \rho_a (\vath)^j \int_0^\infty (\rmvh)^{j+2} \colhla (\rmvh,t) \rmd \rmvh 
      ,\label{dtflweak}
\end{eqnarray}
where the mass density $\rho_a = m_a n_a$. In particular, applying Eqs.\EQ{Mhjl}-\EQo{Rhjl}, we can derive the following relations from the weak form of the FPRS collision spectrum equation\EQ{dtflweak}, reads:
  \begin{eqnarray}
      \pDt \nha \left(t \right) &:=& \frac{\pdt n_a}{n_a} \ = \  \calRh_{0,0},   \label{Rhna} 
      \\
      \pDt \Izha \left(t \right) &:=& \frac{\pdt I_a}{\rho_a \vath} \ = \ \frac{1}{3} \calRh_{1,1} ,  \label{RhIa}
      \\
      \pDt \Kha \left(t \right) &:=& \frac{\pdt K_a}{n_a T_a} \ = \ \frac{3}{2} \calRh_{2,0} ~. \label{RhKa}
  \end{eqnarray}
Applying Eqs.\EQ{Rhna}-\EQo{RhKa}, one can derive the following relation:
  \begin{eqnarray}
      \pDt \Tha \left(t \right) &=& \pDt \Kha \left[2 \Izha  \left(\pDt \Izha - \Izha \Rdtvath  \right) + \Kha \Rdtvath \right] \ \equiv \ 0 ~.  \label{RhTa} 
  \end{eqnarray}
The above equation serves as a convergence criterion for our algorithm in solving the FPRS collision spectral equation\EQ{dtfl}.

Mass, momentum and energy conservation\EQ{Cn}-\EQo{CK} can be reformulated as a function of $\colhlab$\EQgive{collab}, reads:
  \begin{eqnarray}
      \left<1, \colhoab \right>_{\rmvh} &=& \left<1, \colhoba \right>_{\rmvhb} \quad = \quad 0, \label{Cnh}
      \\
      \rho_a \vath \left<\rmvh, \colhIab \right>_{\rmvh} &=& - \rho_b \vbth \left<\rmvhb, \colhIba \right>_{\rmvhb}, \label{CIh}
      \\
      \frac{\rho_a \left(\vath\right)^2}{2}\left<\left(\rmvh \right) ^2, \colhoab \right>_{\rmvh} &=& - \frac{\rho_b \left(\vbth\right)^2}{2} \left<\left(\rmvhb \right)^2, \colhoba \right>_{\rmvhb}, \label{CKh}
  \end{eqnarray}
where $\left<\left(\rmvh\right)^j,g \right>_{\rmvh}$ denotes the integral of $4 \pi (\rmvh)^{j+2} \cdot g$ with respect to $\rmvh$. These conservation constraints are activated when enforcing discrete conservation (details provided in Sec.\SEC{Conservation enforcing}) of the normalized FPRS collision operator \EQgive{FPShdab}. Otherwise, they serve as indicators to evaluate the performance of our algorithm.

\subsection{Speed coordinate}
\label{Speed coordinate}

By selecting sufficient large values of $l_M$ and $L_M$, the spectral equation\EQ{dtfl} can closely approximate the original FPRS collision equation\EQ{FPeqd} with high accuracy. However, due to the time-discrete error and velocity-discrete errors, it is not feasible to directly obtain the proper distribution function for general scenarios which are typically nonlinear multi-scale systems. Nevertheless, phenomena that are more closely related to the conserved moments in both $l$ space and $j$ space, are usually crucial for system evolution and of greater interest to physicists. Therefore, an algorithm that approximates the distribution function based on the convergence of kinetic moments with a specific collection of $(j,l)$ - including mass, momentum and total energy conservation - will be preferred.

Moreover, in the SHE together with KFE framework, the smoothness of $f(\rmv,\mu, t)$ leading to continuous differentiability of all amplitudes, $\fhl(\rmvh,t)$. Additionally, continuity in system evolution over time is also assumed. In contrast to the previous work utilizing the multi-quadric radial basis function\cite{Askari2015MeshlessEquation} , a new function named King (details in Appendix\APP{King Function}) is introduced in speed coordinate. Based on the King function, King method is developed as a moment convergence technology for approximating the desired functions.

\subsubsection{King method}
\label{King method}

The new King function, which is associated with the first class of modified Bessel functions, can be defined as following:
    \begin{eqnarray}
        \Kl \left(\rmvh;\iota,\sigma \right) &=& \frac{(l+1/2)}{\sigma^2 \sqrt{2 \left |\iota \right | \rmvh}} \left(\frac{\left |\iota \right |}{\iota} \right)^l \e^{-\sigma^{-2} \left (\rmvh^2 + \iota^2  \right)} \mathrm{Besseli} \left(\frac{2l+1}{2}, 2 \frac{\left |\iota \right |}{\sigma^2} \rmvh \right) ~.  \label{King}
    \end{eqnarray}
In this context, the independent variable $\rmvh\in [0,\bbR^+]$ and the parameter $l \in [0,\bbN^+]$, representing the order of the King function. The parameters $\iota$ and $\sigma$ are characteristic parameters of the King function, satisfying $\sigma \in \bbR^+$ and $\iota \in \bbR$.

The $l^{th}$-order King function, $\Kl \left(\rmvh;\uzha, \vhath \right)$\EQgive{King} has the same asymptotic behaviour (details in Appendix\APP{King Function}) as the $l^{th}$-order amplitude of the normalized distribution function, $\fhl(\rmvh,t)$\EQgive{fhlfhGQ}. Therefore, we can effectively approximate $\fhl(\rmvh,t)$ using the King function as follows:
    \begin{eqnarray}
        \fhl \left(\rmvh,t \right) &=& \frac{2 \pi} {\pi^{3/2}} \sum{_{s=1}^{\NKa}}  \left [ \nhas \Kl \left(\rmvh;\uzhas,\vhaths \right) \right] ~.  \label{KFE}
    \end{eqnarray}
Name $\uzhas$ as the characteristic group velocity. The number of King function, $\NKa \in \bbN^+$ is a predetermined value at the initial time level but determined by the subsequent optimization scheme (details in Sec.\SEC{Optimization scheme}). Referring Eq.\EQ{KFE} as the King function expansion (KFE), which is a parameter model for the amplitude function. The theoretical convergence of KFE can be established by utilizing Wiener's Tauberian theorem\cite{Wiener1932} , and a heuristic proof is provided in a subsequent work\cite{wang2024Relaxationmodel} posted on arXiv.

Multiplying both sides of Eq.\EQ{KFE} by $4 \pi \rmvh^{j+2} \rmd \rmvh$ and integrating over the semi-infinite interval $\rmvh = \left [0, \infty \right)$ yields the following characteristic parameter equations (CPEs):
      \begin{eqnarray}
            \calMhjl \left(t \right) \ = \
          \left \{
            \begin{aligned}
            & {C_M}_j^l \sum_{s=1}^{\NKa} \nhas (\vhaths)^j \left(\frac{\uzhas}{\vhaths} \right)^l
            \left [1 + \sum_{\gamma=1}^{j/2} C_{j,l}^\gamma \left (\frac{\uzhas}{\vhaths} \right)^{2\gamma} \right], \quad l \in 2 \bbN, \label{MhKj00sM2}  
            \\
            &  {C_M}_j^l \sum_{s=1}^{\NKa} \nhas (\vhaths)^j \left(\frac{\uzhas}{\vhaths} \right)^l 
            \left [1 + \sum_{\gamma=1}^{(j-1)/2} C_{j,l}^\gamma \left (\frac{\uzhas}{\vhaths} \right)^{2\gamma} \right] , \quad l \in 2 \bbN + 1 , \label{MhKj10sM2}  \\
            \end{aligned}  \label{CPEs}
          \right.
      \end{eqnarray}
where $j \in \left\{(l+2j_p-2)|j_p \in \bbN^+ \right\}$. The coefficients $C_{j,l}^\gamma$ and ${C_M}_j^l$ are given by:
      \begin{eqnarray}
          C_{j,l}^\gamma &=& 2^\gamma \frac{(2 l + 1)!!  C_\gamma^{(j-l)/2}}{(2 l + 2 \gamma + 1)!!},
          \quad
          {C_M}_j^l \ = \ \frac{1}{2^{(j-l)/2}}  \frac{(l + j + 1)!!}{(2 l - 1)!!} ~.
      \end{eqnarray}
Symbol $C_\gamma^{(j-l)/2}$ denotes the binomial coefficient. In particular, when all
$\uzhas$ are zero, the CPEs will be reduced to:
      \begin{eqnarray}
            \calMhjl (t) &=&  {C_M}_j^l \sum{_{s=1}^{\NKa}} \nhas (\vhaths)^j, \quad l \in 2 \bbN~. \label{MhKj00sM2KMM}
      \end{eqnarray}
Note that the original definition of normalized kinetic moment\EQ{Mhjl}, $\Mh_{j,l}$, represents the numerical value calculated from the amplitude function before being smoothed by King function in this research. This is accomplished through the utilization of a Romberg integral\cite{Bauer1961AlgorithmIntegration} given by Eq.\EQ{MhjlRI}. While $\calMhjl$ denotes the desired form in physics derived from the KFE model\EQ{KFE}. The relative optimization error is defined as:
\begin{eqnarray}
    \delta \calMhjl (t)  &=& \left |\hat{M}_{j,l} - \calMhjl \right| / \left |\hat{M}_{j,l} \right|, \label{dMhjl}
\end{eqnarray}
which will be utilized to assess the convergence of the optimization process.

Due to the rapid damping of the higher-order harmonics, the first few harmonics in SHE\EQ{fhfhl} contain much of the most important physics for many plasma physical problems of interest. Therefore, we can assume that the characteristic parameters $\nhas$, $\uzhas$ and $\vhaths$ in KFE\EQ{KFE} are independent of $l$. This assumption allows us to approximate the normalized function $\fh \left(\rmvh, \mu,t \right)$ by $\NKa$ Gaussian functions, $\fh \left(\rmvh, \mu,t \right) = \sum_{s=1}^{\NKa}   \frac{1}{\pi ^{3/2}} \frac{\nhas}{\vhaths^3} \e^{-\left (\vh - \uzhas \bfe_z \right )^2}.$ This approximation is reasonable and efficient, particularly for weakly anisotropic plasmas in which each sub-component $\Kl$ is not far from a local equilibrium state (the characteristic group velocity, $|\uzhas| \ll 1$). Therefore, the distribution function is characterized by a total of $3\NKa$ unknown parameters. Given knowledge of any specified $3 N_{k_a}$ kinetic moments $\calMhjl$ with a collection of $(j,l)$, then solving the corresponding well-posed CPEs\EQ{CPEs} can provide us with all the characteristic parameters.
For weakly anisotropic plasmas, the following section provides three specific optimization schemes, L01jd2nh, L01jd2 and L01jd2NK, all of which are convergent for higher-order moments.

\subsubsection{Optimization scheme for characteristic parameters}
\label{Optimization scheme}

In $l$ space, it is preferable to choose the kinetic moments of the first two orders amplitude functions for weakly anisotropic plasmas, due to the simplicity of the corresponding CPEs\EQ{CPEs}. Following the principle of simplicity, we propose the following optimization scheme by selecting the collection of kinetic moments in $j$ space,
      \begin{equation}
          (j,l) \in \left\{(2j_p+l,l)|j_p \in \bbN,0\leq l\leq1,0\leq 2j_p+l<3\NKa\right\}~.
           \label{L01jd2nh}
      \end{equation}
This will be referred as the \textbf{L01jd2nh} scheme. In particular, when $(j,l)=(1,1)$, CPEs represented by Eq.\EQ{MhKj10sM2} yields:
      \begin{eqnarray}
            \sum{_{s=1}^{\NKa}} \nhas \uzhas & =  & \calMhII / 3 \ = \ \uzha ~. \label{uzhanuTs}
      \end{eqnarray}
When $(j,l)=(0,0)$ and $(j,l)=(2,0)$, CPEs gives:
      \begin{eqnarray}
            \sum{_{s=1}^{\NKa}} \nhas & =  & \calMhoo \ = \ \nha, \label{nhanuTs}
            \\
            \sum{_{s=1}^{\NKa}} \nhas \left[3 (\vhaths)^2 /2+ (\uzhas)^2 \right]  &=& \calMh_{2,0} \ = \ {3}/{2} + \left(\uzha \right)^2 ~ . \label{KhanuTs}
      \end{eqnarray}
Eqs.\EQ{uzhanuTs}-\EQo{KhanuTs} serve as conservation constraints and will be utilized as the constraint equations for L01jd2nh, as well as the subsequent L01jd2 scheme. This indicates that when $\NKa\geq 1$, the L01jd2nh scheme can ensure conservation of mass, momentum and energy during the optimization process.

When the parameters $\nhas$ of species $a$ are all nearly constant, a more effective scheme named \textbf{L01jd2}, rather than L01jd2nh, can be utilized. Here, assuming that $\nhas$ are constants and any $2\NKa$ normalized kinetic moments, $\calMhjl$ could be employed in CPEs\EQ{CPEs} to determine the parameters $\uzhas$ and $\vhaths$. In the L01jd2 scheme, $\NKa$ will remain constant and
      \begin{equation}
          (j,l) \in \left\{(2j_p-l,l)|j_p \in \bbN^+,0\leq l\leq1,1\leq j_p<\NKa\right\}~.
           \label{L01jd2}
      \end{equation}
      
The \textbf{L01jd2NK} scheme is obtained by integrating the L01jd2nh scheme and the L01jd2 scheme, using the criterion given by:
      \begin{equation}
          \delta \nhas \ = \ \left|\nhas(\tkI) - \nhas (\tk) \right| / \nhas (\tk), \quad s=1,2,\cdots ,\NKa, \label{dnhaskI}
      \end{equation}
where the superscripts, $k$ and $k+1$, denote the time levels. When the characteristic parameters of species $a$ satisfy $(\NKa)^{-1} \sum_{s=1}^{\NKa} \delta \nhas \le rtol_{n}$, where $rtol_{n}$ is a specified relative tolerance, the L01jd2 scheme will be executed. Otherwise, the L01jd2nh scheme will be performed within L01jd2NK scheme. Unless otherwise stated, the parameter $rtol_{n}=0.1$.  

Especially, when all parameters $\uzhas$ are zero, the L01jd2nh scheme simplifies to:
      \begin{equation}
          (j,l) \in \left\{(2j_p-2,l)|j_p \in \bbN^+,1\leq j_p \le 2\NKa,l=0\right\}
           \label{L01jd2NKKMM}
      \end{equation}
and the L01jd2 scheme becomes:
      \begin{equation}
          (j,l) \in \left\{(2j_p,l)|j_p \in \bbN^+,1\leq j_p \le \NKa,l=0\right\}~.
           \label{L01jd2KMM}
      \end{equation}
Meanwhile, the conservation constraints represented by Eqs.\EQ{uzhanuTs}-\EQo{KhanuTs} are reduced to:
      \begin{equation}
          \calMh_{1,1} = 0, \quad \calMhoo \ = \ 2 \calMh_{2,0} / 3 = 1 ~.
           \label{m}
      \end{equation}
Eq.\EQ{L01jd2NKKMM} and Eq.\EQ{L01jd2KMM} show that if $\NKa$ increases by one, the convergence order $j$ increases by two in the L01jd2 scheme and by four in the L01jd2nh scheme when $\NKa\geq 2$.

After specifying the optimization scheme with a collection of $(j,l)$ and computing the normalized kinetic moments $\Mh_{j,l}$, we can approximate $\calMhjl$ by let $\calMhjl = \Mh_{j,l}$. Hence, the parameters $\nhas$,  $\uzhas$ and $\vhaths$ can be solved from the well-posed CPEs\EQ{CPEs} by utilizing a least squares method\cite{Fong2011LSMR:Problems} (LSM), specifically Levenberg-Marquardt\cite{Wright1985AnSqures,Kanzow2004Levenberg-MarquardtConstraints} method. The main procedure for the King method is outlined in Algorithm \ref{alg: King}. In order to enclose the nonlinear FPRS collision spectral equation\EQ{dtfl}, the approximation of $\fhl$ when $l \ge 2$ can be achieved through the utilization of KFE\EQ{KFE} with established characteristic parameters. The remaining question is how to determine the number of King functions, $\NKa$, at a new time level. This will be provided in the following section.

    \begin{algorithm}[!ht]
    \caption{Employing King method to smooth the normalized amplitude functions of species $a$} \label{alg: King}
      From inputs $l_M(t)$, $\NKa(t)$, $\rmvhA(t)$, collection of $(j,l)$ and $\fhl(\rmvhA(t),t)$
  
      \quad 1 Compute the normalized kinetic moments $\Mh_{j,l}(t)$\EQ{MhjlRItk1}

      \quad 2 Let $\calMhjl := \Mh_{j,l}(t)$
      update $\nhas(t)$, $\uzhas(t)$, $\vhaths(t)$ by solving CPEs\EQ{CPEstk1} 
      
      \quad 3 Evaluate the new number of King functions, $\NKa(t^*)$ according to algorithm \ref{alg: NKa}

      \quad \textbf{If} $\NKa(t^*) \ne \NKa(t)$

          \quad \quad \quad 4 Let $\NKa(t) := \NKa(t^*)$

          \quad \quad \quad 5 Update the collection of $(j,l)$ based on $\NKa(t^*)$

          \quad \quad \quad 6 Solve the CPEs\EQ{CPEstk1} again
          
      \quad \textbf{End}

      \quad 7 Evaluate the convergence of optimization by calculate $\delta \calMhjl$ according to Eq.\EQ{dMhjl}

      \quad 8 Update $l_M(t)$ and the smoothed amplitudes, $\fhl(\rmvhA(t),t)$ according to Eq.\EQ{KFEtk1}

      \quad 9 \textbf{Return} $l_M(t)$, $\NKa(t)$, $\rmvhA (t)$, collection of $(j,l)$ and $\fhl(\rmvhA(t),t)$
    \end{algorithm}

\subsubsection{Update $\NKa$}
\label{Update NK}

The adaptability of $\NKa$ is essential for capturing the nonlinear effects of plasma systems. Before providing the self-adaptive scheme for $\NKa$, we introduce an indistinguishable condition represented by Eq.\EQ{Dvhths} for KFE\EQ{KFE}. If two known groups of characteristic parameters $\left (\iota_1, \sigma_1 \right) $ and $\left (\iota_2, \sigma_2 \right) $ in KFE, each with respective weights $\nhaI$ and $\nhaII$, satisfy
    \begin{eqnarray}
      \left| \frac{\sigma_1}{\sigma_2} - 1 \right| + \left| \frac{\iota_1}{\iota_2} - 1 \right|  & \le & rtol, \label{Dvhths}
    \end{eqnarray}
we claim that the King function $\Kl \left(\rmvh;\iota_1,\sigma_1 \right)$ and $\Kl \left(\rmvh;\iota_2,\sigma_2 \right)$ are identical with parameters $(\iota_0,\sigma_0)$, within the allowed error range. Here, $rtol$ is a specified relative tolerance with a default value of $rtol=10^{-10}$. The weight of $\Kl \left(\rmvh;\iota_0,\sigma_0 \right)$ is determined by $\nhao = \nhaI + \nhaII$.

The number of King functions at the $(k+1)^{th}$ time level will be determined by the given scheme until reaching a specified minimum value, $\NKmin$, or maximum value, $\NKmax$. For weakly anisotropic plasmas, a effective scheme to update $\NKa$ is provided as follows:
    \begin{eqnarray}
      \NKa(\tkI) &=& \NKa(\tk) - d\NKa(\tk),  \label{NK}
    \end{eqnarray}
where
    \begin{eqnarray}
      d\NKa(\tk) &=&
        \begin{aligned}
            \left \{
            \begin{array}{cc}
                0, & Eq.\EQ{Dvhths} == false,  \\
                1, & Eq.\EQ{Dvhths} == true  ~.
            \end{array}
            \right .
        \end{aligned}
         \label{dNK}
    \end{eqnarray}
The value of $\NKa$ at the first time level is set to a specified constant, $\NK$, where  $\NKmin \le \NK \le \NKmax$. The values of $\NKmin=1$ and $\NKmax=3$ have been achieved at present.
    
When $\NKa(\tkI) < \NKa(\tk)$, we can eliminate the higher-order CPEs, which are usually more complex, to obtain an updated well-posed CPEs. The procedures for determining $\NKa(\tkI)$ is outlined in Algorithm \ref{alg: NKa}. However, it should be noted that the method for updating $\NKa$ is not limited to the aforementioned method. A more advanced technique, capable of self-adaptive increase in $\NKa$ with the potential for $\NKmax$ to reach up to 10, is presented in our subsequent work\cite{wang2024Relaxationmodel} and will be further refined in our future research.

    \begin{algorithm}
    \caption{Determine parameter $\NKa$ at the $(k+1)^{th}$ time level.} \label{alg: NKa}
      \textbf{If} optimization scheme == L01jd2

        \quad \quad $\NKa(\tkI) \equiv \NKa (\tko)$
          
      \textbf{Else}
      
      \quad \quad \textbf{If} $k = 0$
      
          \quad \quad \quad \quad $\NKa(\tkI) = \NK$
          
      \quad \quad \textbf{Else}
      
      \quad \quad \quad \quad \textbf{If} Eq.\EQ{Dvhths} == \textbf{True}
          \ $d\NKa(\tk) = 1$
          \ \textbf{Else}
          \ $d\NKa(\tk) = 0$
          \ \textbf{End}
      
           \quad \quad \quad \quad  $\NKa(\tkI) = \min \{\NKmax, \max [\NKmin,\NKa(\tk) - d\NKa(\tk)]\}$ \label{NKatk1}
      
      \quad \quad \textbf{End}
      
      \textbf{End}

    \end{algorithm}

\subsubsection{Discretization of speed coordinate}
\label{Discretization of speed coordinate}

In order to calculate the normalized kinetic moments $\Mh_{j,l}$\EQ{Mhjl} and the Shkarofsky's integrals (give in Eqs.\EQ{IjFL}-\EQo{JjFL}), a set of field nodes in speed coordinate is required. In this research, we employ uniform field nodes in a normalized speed domain of $[0, \rmvh_{M}]$, where $\rmvh_{M}$ denotes the maximum value of the normalized speed. If not specified otherwise, $\rmvh_{M}$ is a constant with a default value of 10. The total number of the field nodes, $N_n=2^{n_2}+1$ where $n_2 \in \bbN^+$. Let $\rmvhA$ represent the $\alpha^{th}$ field node and $[\rmvhA]$ denote the field nodes set where $\alpha = 1,2,\cdots,N_n$. The spacing, $\DvhA$, is then determined as $\DvhA = \rmvh_{M} / 2^{n_2}$. The default value for parameter $n_2$ are set at $n_2=7$, unless otherwise stated. Therefore, the default number of nodes, $N_n$ will be $129$.
  
\begin{figure}[htp]
	\begin{center}
		\includegraphics[width=0.8\linewidth]{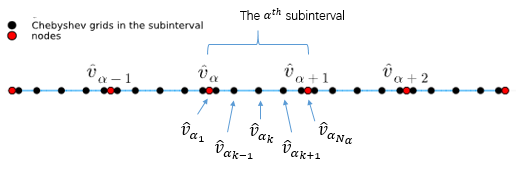}
	\end{center}
	\caption{Discretization of speed coordinate: field nodes and the Chebyshev grids in the subinterval.}
	\label{subinterval}
  \end{figure}

The normalized kinetic moments $\Mh_{j,l}$\EQ{Mhjl} will be calculated using the Romberg integral\cite{Bauer1961AlgorithmIntegration} method. For convenience, the value and its relative error will be expressed as:
  \begin{eqnarray}
      \Mh_{j,l} (t), Error \left(\Mh_{j,l}\right) &=&  \left <([\rmvhA])^{j} | \fhl ([\rmvhA],t) \right > _{R}, \label{MhjlRI}
  \end{eqnarray}
denoting the Romberg integral of function $(\rmvh)^{j+2} \fhl (\rmvh)$ over the set $[\rmvhA]$. Symbol $Error \left(\Mh_{j,l}\right)$ represents the estimated upper bound of the integral error of $\Mh_{j,l}$. The set $[\rmvhA]$ represents the uniform nodes (depicted as red points in Fig.\FIG{subinterval}). Similarly, the computation of Eq.\EQ{Rhjl} will also be performed using Romberg integral:
  \begin{eqnarray}
      \calRhjl (t), Error \left(\calRhjl \right) &=&  \left <([\rmvhA])^{j} | \colhla ([\rmvhA],t) \right > _{R} ~. \label{RhjlRI}
  \end{eqnarray}
The relative errors, $Error \left(\Mh_{j,l}\right)$ and $Error \left(\calRhjl \right)$  can serve as indicators to evaluate the quality of field nodes determined by parameters $(n_2,N_0)$. These indicators could be utilized in constructing a self-adaptive scheme which will be developed in the future.
  
Similar to Eq.\EQ{KFE}, the $l^{th}$-order normalized amplitude function of the background species $b$, $\FhL \left(\rmvhb,t \right)$, can be approximated as:
    \begin{eqnarray}
        \FhL \left(\rmvhb,t \right) &=& \frac{2 \pi} {\pi^{3/2}} \sum{_{r=1}^{\NKb}}  \left [ \nhbr \KL \left(\rmvh;\uzhbr,\vhbthr \right) \right]~.  \label{FhLnuTr}
    \end{eqnarray}
The number of King functions for species $b$, $\NKb \in \bbN^+$, may not be equivalent to the one of species $a$, $\NKa$. The Shkarofsky's integrals, given in Eqs.\EQ{IjFL}-\EQo{JjFL}, involve variable upper/lower bound integration and remain challenging to compute analytically. In this study, they are numerically evaluated using a set of refined grids, as illustrated in Fig.\FIG{subinterval}. 

To create the $\alpha^{th}$ subinterval $[\rmvh_{{\alpha}k}],k=1,2,\cdots,N_{\alpha}$, we add $N_{\alpha}-2$ auxiliary grids (black points) into the $\alpha^{th}$ interval $[\rmvhA,\rmvh_{\alpha+1}]$, as showed in Fig.\FIG{subinterval}. For background species $b$, let $z_{\alpha_k} = \vabth \rmvh_{{\alpha}_k}$ where $z = \vvbth$ and thermal velocity ratio $\vabth=\vath/\vbth$ . The Shkarofsky's integrals can be calculated using parallel Clenshaw-Curtis (CC) quadrature, which is a type of Gauss-Chebyshev quadrature\cite{Press2007NumericalRecipes} . This method can be represented as follows:
  \begin{eqnarray}
      \IjFhL \left(z_{\alpha},t \right) &=& 
        \begin{aligned}
            \left \{
            \begin{array}{cc}
                0, &  \alpha = 1, \\
                \frac{1}{(z_{\alpha})^j} \sum_{s=2}^{\alpha} \left <\left(
                [z_{s_k}] \right)^{j} | \FhL ([z_{s_k}],t) \right > _{CC}, & 2 \le \alpha \le 2^{n_2}+1,
            \end{array}
            \right .
        \end{aligned}
         \label{IjFLMcc}
         \\
      \JjFhL \left(z_{\alpha},t \right) &=& 
        \begin{aligned}
            \left \{
            \begin{array}{cc}
                0, &  \alpha = 2^{n_2}+1, \\
                (z_{\alpha})^j \sum_{s=\alpha}^{2^{n_2}} \left <\left([z_{s_k}] \right)^{-j} | \FhL ([z_{s_k}],t) \right > _{CC}, & 1 \le \alpha \le 2^{n_2}  ~.
            \end{array}
            \right .
        \end{aligned}
         \label{JjFLMcc}
  \end{eqnarray}
Here, the Clenshaw-Curtis quadrature on the $\alpha^{th}$ subinterval can be formulated as:
  \begin{eqnarray}
    \left <([z_{\alpha_k}])^j | g ([z_{\alpha_k}]) \right > _{CC} &=& \sum{_{k = 1}^{N_{\alpha}}} w_k (z_{\alpha_k})^{j+2} g (z_{\alpha_k}) ~. \label{CC}
  \end{eqnarray}
The set $[z_{\alpha_k}]$ are the Gauss collection on the $\alpha^{th}$ quadrature domain where $k=1,2,\cdots,N_{\alpha}$, while $z_{\alpha_k}$ denotes the $k^{th}$ point on the $\alpha^{th}$ subinterval showed in Fig.\FIG{subinterval}. The corresponding integral weight $w_k$ is calculated according to Fornberg's algorithm\cite{Fornberg1998CalculationFormulas} . Function $\FhL\left(\vvbth ,t \right)$ denotes the normalized background distribution function of species $a$. This function can be easily obtained by mapping the analytical function, $\FhL \left(\rmvhb,t \right)$ given by Eq.\EQ{FhLnuTr}, onto the normalized speed coordinate $\vvbth$. This process of mapping can be expressed as:
  \begin{eqnarray}
    \FhL\left(\vvbth ,t \right) &=& \FhL(\vabth \rmvh,t) ~. \label{mapping}
  \end{eqnarray} 
  
In this paper, the number of auxiliary Chebyshev grids in subintervals, $N_{\alpha}$, is fixed at a given constant $N_0$ and the auxiliary grids (depicted as black points in Fig.\FIG{subinterval}) are only utilized in the step for calculating the Shkarofsky's integrals. Consequently, the maximum number of grids in speed coordinate is determined by $N_{\rmv}=(N_0-1)(N_n-2)+N_0$ and $N_{\rmv 2}=\lMI\times N_{\rmv}$ for the total velocity space. Given the default value of $N_0=7$ unless otherwise specified, we obtain $N_{\rmv}=769$. It should be noted that $N_{\rmv}$ and $N_{\rmv 2}$ are primarily utilized in the Rosenbluth potentials step (Sec.\SEC{Rosenbluth potentials}), while the maximum number of nodes for other steps is determined by number of field nodes, $N_n$. 

It is particularly noteworthy that for moderately anisotropic plasma systems in which some characteristic group velocity $|\uzhas|$ near one, including the so-called subsonic and low supersonic regions, the characteristic parameters in KFE\EQ{KFE} typically are dependent on $l$. The framework of SHE together with KFE remains suitable for these general scenarios. However, the optimization scheme need to be extended as a self-adaptive version, incorporating a comprehensive collection of $(j,l)$, a clearer scheme for updating $\NKa$ and self-adaptive field nodes determined by parameters $(n_2, N_0)$. This is beyond the scope of this paper and will be addressed in future research. Additionally , it is worth mentioning that the KFE serves as a smoothing step, combined with the implicit time discretization given in Sec.\SEC{Time integration}, enabling our algorithm to surpass the classical CFL condition limit. From now on, the subscript and superscript "k" will represent the time level.

\section{Implicit temporal discretization and nonlinear conservation constraints}
\label{Time integration}

Given the specified discretization of speed coordinate determined by the value of $\rmvh_M$ and parameters $(n_2,N_0)$, the semi-discrete FPRS collision spectrum equations\EQ{dtfl} at the $\alpha^{th}$ node can be formulated as:
  \begin{eqnarray}
      \ddt \fl \left(\rmvhA,t \right) &=& \colla \left(\rmvhA,t \right) , \quad 0 \le l \le l_M ~. \label{dtfla}
  \end{eqnarray}
where
  \begin{eqnarray}
      \fl \left(\rmvhA,t \right) &=& \navth \fhl \left(\rmvhA,t \right), \label{fla}
      \\
      \colla \left(\rmvhA,t \right) &=& \navth \sum{_{b=1}^{N_s}} \nbvth \Gab \colhlab, \label{Cla}
      \\
       \colhlab \left(\rmvhA, t \right) &=& \sum{_{\beta=1}^{\lMI}} \wmub \Pl (\mub) \colhab \left (\rmvhA, \mub, t \right)
       ~. \label{collaba}
  \end{eqnarray}
Here, function $\colhab$ is a nonlinear function of $\fhl$ and $\FhL$ and satisfies the form presented in Appendix\APP{Normalized FPRS collision operator}. Number density remains constant theoretically and thermal velocity, $\vath$, is a function of $\rho_a$, $I_a$\EQ{Iha} and $K_a$\EQ{Kha}, reads:
\begin{eqnarray}
    \vath &=& \sqrt{\frac{2}{3} \left(\frac{2 K_a}{\rho_a} - \left(\frac{I_a}{\rho_a} \right)^2 \right)} ~.\label{vath}
\end{eqnarray}

The final semi-discrete scheme for the nonlinear FPRS collision equations is derived by incorporating the aforementioned equations, including the FPRS collision spectrum equations\EQ{dtfla}, the FPRS collision operators represented by Eqs.\EQ{Cla}-\EQo{collaba} and\EQ{FPShdabSi}. The King function expansion\EQ{KFE} acts as a smoothing step, while the characteristic parameter equations\EQ{CPEs}, convergence criterion\EQ{RhTa} and conservation equations\EQ{Cnh}-\EQo{CKh} serve as constraints for FPRS collision spectrum equations\EQ{dtfla}.

\subsection{Implicit temporal discretization}
\label{Implicit temporal discretization}

\subsubsection{Time block technique}
\label{Time block technique}

The discrete form of FPRS collision spectrum equation\EQ{dtfla} utilizing explicit Euler scheme\cite{Rackauckas2017DifferentialEquations.jlaJulia} can be formulated as follows:
  \begin{eqnarray}
      \fl \left(\rmvhA^{k},\tkI \right) &=& \fl \left(\rmvhA^{k},\tk \right) + \Dtk \colla \left(\rmvhA^{k},\tk \right) , \label{fltk1ExEuler}
  \end{eqnarray}
where $\Delta \tk=\tkI-\tk$ represents the current timestep size. Similarly, 
  \begin{eqnarray}
      X_a \left(\tkI \right) &=& X_a \left(\tk \right) + \Dtk 
      \partial_t
      X_a \left(\tk \right), \quad X=n,I,K  ~. \label{nIKatk1ExEuler}
  \end{eqnarray}
 After computing ${\vathkI}$\EQ{vath} and $\fhl \left(\rmvhA^{k},\tkI \right)$\EQ{fla}, applying the King method provided by Algorithm \ref{alg: King} yields the analytical formula of $\fl \left(\rmvhA^{k},\tkI \right)$. Notes that $\fl \left(\rmvhA^{k},\tkI \right)$ in Eq.\EQ{fltk1ExEuler} is evaluated on the normalized speed coordinate at the $k^{th}$ time level,  $\rmvhA^{k}=\rmvA^{k}/\vathk $.
 
 Regard the interval $[\tk,\tkI]$ as the $k^{th}$ time block. The time block technique (TBT) suggests that within the $k^{th}$ time block, $\rmvA(t)$ remains invariant and will be updated at the $(k+1)^{th}$ time level within the $(k+1)^{th}$ time block. The default update will keep $\rmvhA^{k+1}$ equivalent to $\rmvhA^{k}$, determined by $\rmvh_M$ and $(n_2,N_0)$. Defining $\rmvhA^{k,k+1}=\rmvA^{k}/\vathkI$ gives
  \begin{eqnarray}
      \rmvhA^{k,k+1} &=& \rmvhA^{k} {\vathk}/{\vathkI},   \label{vhkk1}
  \end{eqnarray}
Hence, the amplitude function on $\alpha^{th}$ node at $(k+1)^{th}$ time level within $k^{th}$ time block is:
  \begin{eqnarray}
      \fl \left(\rmvhA^{k,k+1},\tkI \right) &=& \fl \left(\rmvhA^{k} {\vathk}/{\vathkI},\tkI \right)  ~. \label{flk1map}
  \end{eqnarray}
The successful application of TBT relies on the analytical KFE\EQ{KFE} in King method.

The explicit Euler scheme is a first-order precision algorithm, which can be employed as the prediction step (first stage) of an implicit algorithm, such as the trapezoidal scheme. Moreover, TBT is also effective in Trapezoidal scheme, presented in Algorithm \ref{alg: Trapez}. The main procedure for solving FPRS collision spectrum equation within $k^{th}$ time block utilizing the explicit Euler scheme is outlined in the following pseudo-code (Algorithm \ref{alg: ExEuler}).

    \begin{algorithm}[!ht]
    \caption{Solving FPRS collision spectrum equation by utilizing explicit Euler scheme} \label{alg: ExEuler}
      From inputs $\Dtk$, $l_M (\tk)$, $\NKa(\tk)$, $\rmvhA^{k}$,collection of $(j,l)$ and $\fl(\rmvhA^{k},\tk)$ for all species

      \quad 1 Calculate $\ddt \nak$, $\ddt \Izak$, $\ddt \Kak$, $\rmvhA^{k}$ and $\colla (\rmvhA^{k},\tk)$ according to Algorithm \ref{alg: Ca}

      \quad 2 Evaluate $\fl (\rmvhA^{k},\tkI)$ according to Eq.\EQ{fltk1ExEuler}
      
      \quad 3 Calculate ${\nakI}$, ${\IzakI}$, ${\KakI}$ according to Eq.\EQ{nIKatk1ExEuler} 

      \quad 4 Update ${\vathkI}$ according to Eq.\EQ{vath} and $\fhl(\rmvhA^{k},\tkI)$ according to Eq.\EQ{fla} 

      \quad 5 Smooth $\fhl(\rmvhA^{k},\tkI)$ utilizing King method presented by Algorithm \ref{alg: King}
      
      \quad 6 Update $\rmvhA^{k,k+1}$ according to Eq.\EQ{vhkk1} and $\fl(\rmvhA^{k,k+1}{},\tkI)$ according to Eq.\EQ{flk1map}

      \quad 7 \textbf{Return} $\ddt \nak$, $\ddt \Izak$, $\ddt \Kak$, $\colla (\rmvhA^{k},\tk)$, ${\vathkI}$, $\rmvhA^{k,k+1}$ and $\fl(\rmvhA^{k,k+1},\tkI)$

    \end{algorithm}

\subsubsection{Implicit iteration}
\label{Implicit iteration}

For equation $\partial_t A = g(A)$, the discretization employing trapezoidal\cite{Rackauckas2017DifferentialEquations.jlaJulia} scheme can be formulated as follows:
\begin{eqnarray}
    A^{k+1} &=& A^{k} + \Dtk \left[c_{k} g \left(A^k \right) + c_{k+1} g \left({A^{k+1}}^* \right) \right], \label{Ak1}
\end{eqnarray}
where the coefficients $c_k$ and $c_{k+1}$ represent the temporal weights for Range-Kutta method. In the trapezoidal scheme, both of them are equal to $1/2$. During the implicit iteration to optimize the convergence of $A$ at $(k+1)^{th}$ time level, ${A^{k+1}}$ denotes the value of $A$ at the $i^{th}$ stage of $(k+1)^{th}$ time level, while ${A^{k+1}}^*$ represents the value of $A$ at the ${(i-1)}^{th}$ stage of $(k+1)^{th}$ time level, where $i=2,3,\cdots,N_{in}$. The maximum number of implicit iteration within each time block is denoted as $N_{in}$, typically set to a default value of 10.

By employing the trapezoidal scheme for time integration, the discrete FPRS collision spectral equation\EQ{dtfla} at $\alpha^{th}$ node and within $k^{th}$ time block can be formulated as:
  \begin{eqnarray}
      \fl \left(\rmvhA^{k,k+1},\tkI \right) &=& \fl \left(\rmvhA^{k},\tk \right)
      + 
      \Dtk \left [c_{k} \colla \left(\rmvhA^{k},\tk \right) + c_{k+1} \colla \left(\rmvhA^{k,k+1},\tkI \right) \right]~. \label{fltk1Trapez}
  \end{eqnarray}
The pseudo-code for solving Eq.\EQ{fltk1Trapez} is provided by Algorithm \ref{alg: Trapez}. The $l^{th}$-order amplitude of FPRS collision operator at $\alpha^{th}$ node and $(k+1)^{th}$ time level within $k^{th}$ time block are:
  \begin{eqnarray}
      \colla \left(\rmvhA^{k,k+1},\tkI \right) &=& \navthkI \sum{_{b=1}^{N_s}} \nbvthkI \Gab^{k+1} \colhlab \left(\rmvhA^{k,k+1},\tkI \right) , \label{FPShldatk1}
  \end{eqnarray}
where
  \begin{eqnarray}
      \colhlab \left(\rmvhA^{k,k+1},\tkI \right) &=& \sum{_{\beta=1}^{\lMI}} \wmub^{k+1} \Pl (\mub^{k+1}) \colhab \left (\vh^{k,k+1}, \tkI \right) ~. \label{FPShldabtk1}
  \end{eqnarray}
Here $\vh^{k,k+1} = \v^{k} / \vathkI$. The normalized FPRS collision operator, $\colhab$, satisfies Eq.\EQ{FPShdabSi}, which is a nonlinear model and can be rewritten as:
  \begin{eqnarray}
  \begin{aligned}
      \colhab & \left(\vh^{k,k+1}, \tkI \right) = \ \colhab \left [
      \fh_0 \left(\rmvhA^{k,k+1},\tkI \right), \cdots, \fhl \left(\rmvhA^{k,k+1},\tkI \right), \cdots,\fh_{l_{max}} \left(\rmvhA^{k,k+1},\tkI \right), 
      \right. \\ &  \quad \quad \quad
      \left.
      \Fh_0 \left(\vabth^{k+1} \rmvhA^{k,k+1},\tkI \right), \cdots,\FhL \left(\vabth^{k+1} \rmvhA^{k,k+1},\tkI \right), \cdots \Fh_{l_{max}} \left(\vabth^{k+1} \rmvhA^{k,k+1},\tkI \right) \right]~. 
  \label{FPShdabSitk1} 
  \end{aligned}
  \end{eqnarray}
  
Note that $\NKa$, $\NKb$, $l_{max}$ are also recalculated at each time level. The scheme to update $\NKa$ is provided in Sec.\SEC{Update NK}. The main procedure for computing normalized FPRS collision operator between species $a$ and species $b$ at $(k+1)^{th}$ time level within $k^{th}$ time block is outlined in Algorithm \ref{alg: Chab}, while procedure for calculating the FPRS collision operator of species $a$ at the $(k+1)^{th}$ time level within $k^{th}$ time block is outlined in Algorithm \ref{alg: Ca}.

    \begin{algorithm}[!ht]
    \caption{Calculating the normalized FPRS collision operator between species $a$ and species $b$ at $(k+1)^{th}$ time level within $k^{th}$ time block} \label{alg: Chab}
      
      From inputs $n_a^{k+1}$, $\vath^{k+1}$, $l_M (\tkI)$, $\NKa(\tkI)$,$\rmvhA^{k,k+1}$ and $\fhl(\rmvhA^{k,k+1},\tkI)$ for all species 
      
      \quad 1 Evaluate $\ddrmvh \fhl(\rmvhA^{k,k+1},\tkI)$, $\dddrmvh \fhl(\rmvhA^{k,k+1},\tkI)$ according to KFE\EQ{KFEtk1}
      
      \quad 2 Update $\FhL(\vabth^{k+1} \rmvhA^{k,k+1} ,\tkI)$ according to 
      the process of mapping\EQ{mapping}

      \quad 3 Calculate Shkarofsky's integrals according to Eqs.\EQ{IjFLMcc}-\EQo{JjFLMcc}

      \quad 4 Compute the amplitudes of Rosenbluth potentials according to Eqs.\EQ{HhL}-\EQo{GhL}

      \quad 5 Update derivatives of Rosenbluth potential's amplitudes according to Eqs.\EQ{dvHhL}-\EQo{ddvGhL}

      \quad 6 Calculate normalized FPRS collision operator, $\colhab (\rmvhA^{k,k+1},\tkI)$ according to Eq.\EQ{FPShdabSitk1}

      \quad 7 Update the amplitudes $\colhlab (\rmvhA^{k,k+1},\tkI)$ according to Eq.\EQ{FPShldabtk1}
    \end{algorithm}

Upon discretization, implicit methods for the FPRS collision equation lead to a complex system of nonlinear algebraic equations, necessitating an effective nonlinear solver strategy for its solution. In this study, we rely on the King method for this task. Under the assumption that the characteristic parameters are independent of $l$, the $l^{th}$-order amplitude function of species $a$ at the $\alpha^{th}$ node and ${(k+1)}^{th}$ time level within $k^{th}$ time block will be smoothed by KFE\EQ{KFE}, reads:
\begin{eqnarray}
    \fhl (\rmvhA^{k,k+1},\tkI) &=&  \frac{2 \pi} {\pi^{3/2}} \sum{_{s=1}^{\NKa}}  \left [ \nhas^{k+1} \Kl \left(\rmvhA^{k,k+1};\uzhas^{k+1},\vhaths^{k+1} \right) \right], \quad 0 \le l \le l_{max} ~.  \label{KFEtk1}
\end{eqnarray}
The characteristic parameters are determined by solving the CPEs\EQ{CPEs}, at ${(k+1)}^{th}$ time level within $k^{th}$ time block, as given by:
      \begin{eqnarray}
            \calMhjl \left(\tkI \right) \ = \
          \left \{
            \begin{aligned}
            & {C_M}_j^l \sum_{s=1}^{\NKa} \nhas^{k+1} \left(\vhaths^{k+1}\right)^j \left(\frac{\uzhas^{k+1}}{\vhaths^{k+1}} \right)^l
            \left [1 + \sum_{i=1}^{j/2} C_{j,l}^i \left (\frac{\uzhas^{k+1}}{\vhaths^{k+1}} \right)^{2i} \right], \quad l =0, \label{MhKj00sM2tk1}  
            \\
            &  {C_M}_j^l \sum_{s=1}^{\NKa} \nhas^{k+1} \left(\vhaths^{k+1}\right)^j \left(\frac{\uzhas^{k+1}}{\vhaths^{k+1}} \right)^l 
            \left [1 + \sum_{i=1}^{(j-1)/2} C_{j,l}^i \left (\frac{\uzhas^{k+1}}{\vhaths^{k+1}} \right)^{2i} \right] , \quad l =1 , \label{MhKj10sM2tk1}  \\
            \end{aligned}  \label{CPEstk1}
          \right.
      \end{eqnarray}
where $j$ satisfies the collection of $(j,l)$ provided in Sec.\SEC{Optimization scheme}. Here, the normalized kinetic moment calculated using Romberg integral\EQ{MhjlRI} at the ${(k+1)}^{th}$ time level within $k^{th}$ time block can be formulated as:
  \begin{eqnarray}
      \Mh_{j,l} (\tkI), Error \left(\Mh_{j,l}(\tkI) \right) &=&  \left <\left([\rmvhA^{k,k+1}] \right)^{j} | \fhl ([\rmvhA^{k,k+1}],\tkI) \right > _{R} ~. \label{MhjlRItk1}
  \end{eqnarray}

    \begin{algorithm}[!ht]
    \caption{Calculating the FPRS collision operator of species $a$ at the $(k+1)^{th}$ time level within $k^{th}$ time block} \label{alg: Ca}
      From inputs $\Dtk$, $l_M (\tkI)$, $\NKa(\tkI)$,$\rmvhA^{k,k+1}$, collection of $(j,l)$ and $\fl(\rmvhA^{k,k+1},\tkI)$

      1 Update $\nakI$, $\IzakI$, $\KakI$  
      
      2 Update $\fhl(\rmvhA^{k,k+1},\tkI)$ according to Eq.\EQ{fla} for all species

      \textbf{For} $a=1,2,\cdots,N_s$
      
      \textbf{For} $b=a,2,\cdots,N_s$

      \quad 3 Evaluate  $\colhlab (\rmvhA^{k,k+1},\tkI)$ according to Algorithm \ref{alg: Chab}

      \quad \textbf{If} Conservation enforcing $==$ \textbf{True}

         \quad \quad 4A Update $\ddt \nakI$, $\ddt \IzakI$, $\ddt \KakI$ according to conservation enforcing algorithm \ref{alg: Conservation enforcing}

      \quad \textbf{Else}

         \quad \quad 4B Update $\ddt \nakI$, $\ddt \IzakI$, $\ddt \KakI$ according to Eqs.\EQ{Rhna}-\EQo{RhKa}
      
      \quad \textbf{End}
      
      \textbf{End} (for species $b$)
      
      \textbf{End} (for species $a$)

      5 Calculate total derivatives of conserved moments respective to $t$: $\ddt \nakI$, $\ddt \IzakI$, $\ddt \KakI$

      6 Compute $\colla (\rmvhA^{k,k+1},\tkI)$ according to Eq.\EQ{FPShldatk1}

      7 \textbf{Return} $\ddt \nakI$, $\ddt \IzakI$, $\ddt \KakI$, $\rmvhA^{k,k+1}$ and $\colla (\rmvhA^{k,k+1},\tkI)$
      
    \end{algorithm}

To properly address those nonlinear algebraic equations described by Eqs.\EQ{fltk1Trapez}-\EQo{MhjlRItk1}, it is imperative to impose constraints for ensuring conservation properties. By employing the Romberg integration with appropriate field nodes, we are able to achieve high-precision values for moments, including number density, momentum and energy. Building upon these high-precision values, according to manifold theory, we rely on backward error analysis to guarantee exact conservation. This is accomplished in Sec.\SEC{Conservation enforcing}.

\subsection{Conservation enforcing}
\label{Conservation enforcing}

The accuracy of integrals\EQ{Rhna}-\EQo{RhKa} can be guaranteed for at least one during two-species collision processes (details in Sec.\SEC{Two-species thermal equilibration}). The convergence criterion for conservation with high accuracy in two-species collision processes at the ${(k+1)}^{th}$ time level within $k^{th}$ time block can be expressed as:
  \begin{eqnarray}
       \min(\pDt \hat{C}_a^{k+1}, \pDt \hat{C}_b^{k+1}) &\ll& 1, \label{Ccond}
  \end{eqnarray}
where
  \begin{eqnarray}
       \pDt \hat{C}_a^{k+1} &=& \left|Error(\pDt \nh_a^{k+1}) \right| + \left|Error(\pDt \hat{I}_a^{k+1}) \right| + \left|Error(\pDt \Kh_a^{k+1}) \right| \label{CnIKh}
  \end{eqnarray}
and similar for $\pDt \hat{C}_b^{k+1}$. 
The aforementioned criterion will also be employed to evaluate the quality of field nodes as described in Sec.\SEC{Discretization of speed coordinate}. Note that this current research will not utilize a self-adaptive scheme for $(n_2,N_0)$.

From manifold theory\cite{HairerE.2006GeometricIntegration} , post-step projection onto manifolds maintains a consistent convergence rate, and conservation properties can be preserved as long as the local solution errors remain sufficiently small. Consequently, incorporating a conservation strategy into our algorithm becomes feasible. This strategy enforces discrete conservation equations\EQ{Cnh}-\EQo{CKh}, by utilizing the more precise integrals, such as those for species $b$ described by Eqs.\EQ{Rhna}-\EQo{RhKa}, to provide more accurate representations during two-species collision processes. The convergence of this conservation strategy will occur when the criterion given by Eq.\EQ{Ccond} is satisfied. 

According to the conservation constraints represented by Eqs.\EQ{CIh}-\EQo{CKh}, by applying relations\EQ{Rhna}-\EQo{RhKa}, the rates of momentum and energy change of species $a$ with respect to time at the $(k+1)^{th}$ time level can be theoretically expressed as:
  \begin{eqnarray}
      \ddt \IzakI &=& - \ddt \IzbkI ,  \label{dtIaIbk}
      \quad
      \ddt \KakI \ = \  - \ddt \KbkI ~. \label{dtKaKbk}
  \end{eqnarray}
When enforcing conservation\EQ{Rhna}, all rates of change in number density will be zero, yields:
\begin{eqnarray}
    \nakI &=& \nak. \label{nak1}
\end{eqnarray}

However, it is not feasible to achieve an exact numerical realization of Eq.\EQ{dtIaIbk} using a general integral scheme, such as Romberg integral. When all the local errors , represented by $Error \left(\pDt \hat{I}^{k+1} \right)$ and $Error \left(\pDt \Kh^{k+1} \right)$ are sufficiently small through the utilization of appropriate field nodes determined by $(n_2,N_0)$, we can enforce Eq.\EQ{dtIaIbk} by selecting the more precision one between species $a$ and species $b$. For example, if the local errors of species $b$ are small,  $\pDt \hat{C}_a^{k+1} \ge \pDt \hat{C}_b^{k+1}$, then according to manifold theory, let
  \begin{eqnarray}
      \ddt \IzakI &:=& - \ddt \IzbkI,  \label{dtIaIbk2}
      \quad
      \ddt \KakI \ := \ - \ddt \KbkI \label{dtKaKbk2}
  \end{eqnarray}
at each stage of the $(k+1)^{th}$ time level within $k^{th}$ time block.

In the trapezoidal method, the values of momentum and energy at the $(k+1)^{th}$ time level within $k^{th}$ time block are calculated through an implicit iteration, for species $a$, can be expressed in the following form:
\begin{eqnarray}
    \IzakI = \Izak + \frac{1}{2} \Dtk \left(\ddt {\Izak} +  \ddt {\IzakI}^* \right), \label{Izak1}
    \quad
    \KakI  =  \Kak + \frac{1}{2} \Dtk \left(\ddt {\Kak} + \ddt {\KakI}^* \right) ~.\label{Kak1}
\end{eqnarray}
The implicit iteration at each time level will be terminated when the thermal velocity of all plasma species at the $i^{th}$ stage of ${(k+1)}^{th}$ time level within $k^{th}$ time block satisfies the following condition, for species $a$ reads:
      \begin{eqnarray}
            \left|\frac{\vath(t_{{k+1}_i})}{\vath(t_{{k+1}_{i-1}})} - 1 \right| \le 10^{-6}, \quad i \ge 2~. \label{Cconvergent}
      \end{eqnarray}
      
Consequently, the average velocity and thermal velocity of species $a$ at the $(k+1)^{th}$ time level within $k^{th}$ time block can be expressed as:
\begin{eqnarray}
    \uzakI &=& \frac{\IzakI}{\rhoakI}, \label{uak1}
    \quad
    \vathkI \ = \ \sqrt{\frac{2}{3} \left[\frac{2 \KakI}{\rhoakI} - \left(\uzakI \right)^2 \right]} ~.\label{vathk1}
\end{eqnarray}
Thus, the normalized average velocity of species $a$ is calculated as $\uzhakI=\uzakI/\vathkI$. By applying Eqs.\EQ{uzhanuTs}-\EQo{KhanuTs}, we obtain:
      \begin{eqnarray}
            \calMhoo (\tkI) & \equiv & 1, \label{nhak1nuTs}
            \\
            \calMhII (\tkI) &=& 3 \uzhakI, \label{uzhak1nuTs}
            \\
            \calMh{_{2,0}} (\tkI)  &=& \frac{3}{2} + \left(\uzhakI \right)^2 ~ . \label{Khak1nuTs}
      \end{eqnarray}
 The conservation enforcing algorithm is presented in the following pseudo-code 
 (\ref{alg:  Conservation enforcing}).

    \begin{algorithm}[!ht]
    \caption{Enforcing conservation during two-species collision processes at each stage of $(k+1)^{th}$ time level within $k^{th}$ time block.} \label{alg:  Conservation enforcing}
      From inputs $\colhlab(\rmvhA^{k,k+1}, \tkI)$ of species $a$ and species $b$:

      \quad 1 Evaluate $\calRhjl(\tkI)$ and $Error \left(\calRhjl(\tkI) \right)$ by calculating the Romberg integral\EQ{RhjlRI}
      
      \quad 2 Update $\ddt \nkI$, $\ddt \IzkI$, $\ddt \KkI$ for species $a$ and species $b$ according to Eqs.\EQ{Rhna}-\EQo{RhKa}
      
      \quad 3 Evaluate whether $\min(\pDt \hat{C}_a^{k+1}, \pDt \hat{C}_b^{k+1}) \ll 1$\EQgive{Ccond}
      
      \quad \textbf{If} $\pDt \hat{C}_a^{k+1} \ge \pDt \hat{C}_b^{k+1}$
      
           \quad \quad 4A $\ddt \IzakI := - \ddt \IzbkI$, $\quad \ddt \KakI := - \ddt \KbkI$
           
      \quad \textbf{Else}
      
           \quad \quad 4B $\ddt \IzbkI := - \ddt \IzakI$, $\quad \ddt \KbkI := - \ddt \KakI$

      \quad \textbf{End}
      
      \quad 5 Let $\ddt \nkI := 0$ for all species

      \quad 6 \textbf{Return} $\ddt \nkI$, $\ddt \IzkI$, $\ddt \KkI$ for all species.

    \end{algorithm}

 By implementing the conservation enforcing scheme mentioned above, the accuracy of the conservation\EQ{Cnh}-\EQo{CKh} are determined by the precision of the more accurate species rather than the less accurate one. It is worth noting that small local solution errors of $\pDt \nh$\EQ{Rhna}, $\pDt \Ih$\EQ{RhIa} and $\pDt \Kh$\EQ{RhKa} for at least one species during two-species Coulomb collision process, represented by Eq.\EQ{Ccond}, are a necessary condition for convergence. This condition is verified at each stage and every time level in our algorithm. Finally, the main procedure for solving FPRS collision spectrum equation by utilizing trapezoidal scheme at the $(k+1)^{th}$ time level within $k^{th}$ time block is outlined in the
 Algorithm \ref{alg: Trapez}.

    \begin{algorithm}[!ht]
    \caption{Solving FPRS collision spectrum equation by utilizing trapezoidal scheme at the $(k+1)^{th}$ time level within $k^{th}$ time block} \label{alg: Trapez}
      From inputs $\Dtk$, $\rmvhA^{k}$, $l_M (\tk)$, $\NKa(\tk)$, collection of $(j,l)$ and $\fl(\rmvhA^{k},\tk)$ for all species

      Initial $N_{iter}=1$ and denote $X = n, I, K$
      
      1 Compute $\ddt X_a^{k}$, $\colla (\rmvhA^{k},\tk)$, ${\vathkI}{^*}$,$\rmvhA^{k,k+1}{^*}$ and $\fl(\rmvhA^{k,k+1}{^*},\tkI)$ according to Algorithm \ref{alg: ExEuler}

      \textbf{For} $N_{iter} = 2,3,\cdots,N_{in}$
      
      \quad 2 Update $\ddt X_a^{k+1}{^*}$,$\rmvhA^{k,k+1}{^*}$ and $\colla (\rmvhA^{k,k+1}{^*},\tkI)$ according to Algorithm \ref{alg: Ca}
      
      \quad 3 Let $\colla (\rmvhA^{k,k+1}{^*},\tkI) = \frac{1}{2} [\colla (\rmvhA^{k},\tk) + \colla (\rmvhA^{k,k+1}{^*},\tkI)]$
      
      \quad 4 Let $\ddt X_a^{k+1} = \frac{1}{2} [\ddt X_a^{k} + \ddt X_a^{k+1}{^*}]$

      \quad 5 Evaluate $\fl (\rmvhA^{k,k+1}{^*},\tkI)$ according to Eq.\EQ{fltk1ExEuler}
      
      \quad 6 Calculate $X_a^{k+1}$ according to Eq.\EQ{nIKatk1ExEuler}

      \quad 7 Compute ${\vathkI}$ according to Eq.\EQ{vathk1}

      \quad 8 Calculate $\fhl(\rmvhA^{k,k+1}{^*},\tkI)$ according to Eq.\EQ{fla}

      \quad 9 Smooth $\fhl(\rmvhA^{k,k+1}{^*},\tkI)$ utilizing King method according to Algorithm \ref{alg: King}

      \quad 10 Update $\rmvhA^{k,k+1} = \rmvhA^{k,k+1}{^*} {\vathkI{^*}}/{\vathkI}$ and $\fl(\rmvhA^{k,k+1},\tkI)$ according to Eq.\EQ{nIKatk1ExEuler} 

      \quad 11 \textbf{If} Eq.\EQ{Cconvergent} == \textbf{True}
                  \textbf{Break}
               \textbf{Else} 
                   $\vathkI{^*} = {\vathkI}$ and $\fhl(\rmvhA^{k,k+1}{^*},\tkI) =\fhl(\rmvhA^{k,k+1},\tkI)$
               \textbf{End}

      \textbf{End}
  
      12 Update $\rmvhA^{k+1}$  determined by $\rmvh_M$, $(n_2,N_0)$ and $\fl(\rmvhA^{k,k+1}{^*},\tkI)$ for $(k+1)^{th}$ time block
      
    \end{algorithm}

\subsection{Timestep}
\label{Timestep}

The Coulomb collision process encompasses multiple dynamical times-scales (such as inter-species time-scale, self-collision time-scale, relaxation time-scale of conserved moments, et al.), making it stiff. In this paper, a timestep of $\Delta t = 2^{-5}$ (unless otherwise stated) is utilized for fixed timestep cases. Self-adaptive timestep will be employed (unless otherwise stated) to improve the algorithm performances, which is determined by the following algorithm:
  \begin{eqnarray}
      \DtkI &=& \min \left(ratio_{dt_k} \times \Dtk, ratio_{M_j} \times \left|\frac{1}{y^k} \frac{\partial {y^k}}{\partial t} \right| \right). \label{dtk1_y}
  \end{eqnarray}
Here, the subscripts $k$ denotes the time level and $y$ represents momentum $I$ and total energy $K$ for all species. $\Delta_{t_k}=1$ at the initial time level. Parameters $ratio_{dt_k}$ and $ratio_{M_j}$ are given constants, with default values of $ratio_{dt_k}=1.1$ and $ratio_{M_j}=0.01$ in this paper. For cases utilizing self-adaptive timestep, nearly all timesteps satisfy $10^{-3} \le \Delta t \le 0.1$. 

As a contrast, the explicit timestep size determined by the CFL condition\cite{Taitano2015AEquation} is calculated as follows:
  \begin{eqnarray}
      \DtExp &=& R_{CFL} \times \underset{l = 0,1,\cdots,l_M \atop a,b=1,2,\cdots,N_s}{\min} \left \{\frac{\Dvh}{A_l{_{ab}}}, \frac{(\Dvh)^2}{D_l{_{ab}}} \right \}. \label{dtCFL}
  \end{eqnarray}
In this equation, the parameter $R_{CFL}=0.1$ is utilized in explicit method to ensure long-term stability. Additionally, $A_l{_{ab}}$ and $D_l{_{ab}}$ represent the transport coefficients in the FPRS collision operator\EQgive{FPRS}, reads:
  \begin{eqnarray}
      A_l{_{ab}} &=& \navth \nbvth \times \CHh \Gab \ddscrvh \HhL, \label{Al0ab}
      \quad
      D_l{_{ab}} \ = \ \navth \nbvth \times \CGh \Gab \dddscrvh \GhL ~. \label{Dl0ab}
  \end{eqnarray}
  This explicit timestep size\EQ{dtCFL} will not be directly utilized in our algorithm, but rather serve as a reference for our timestep represented by Eq.\EQ{dtk1_y}.

\section{Numerical results}
\label{Numerical results}

In order to demonstrate the convergence and effectiveness of our method for solving the FPRS collision equation\EQ{FPeqd}, we will assess the performance of our algorithm with various examples of different degrees of complexity. In the benchmarks conducted in this session, the initial distribution functions for particles at $t = 0$ are drifting Maxwellian distributions with a specified number density $n_a$, average velocity $u_a$ and temperature $T_a$, which can be written as:
\begin{eqnarray}
    f(\vh,t)  &=& \navth \fh(\vh,t) \ = \ \frac{1}{\pi ^{3/2}} \navth  \sum{_{s=1}^{\NKo}}  \left [\frac{\nhas} {(\vhaths)^3} \e^{-\left (\vh - \uzhas \bfe_z \right )^2} \right]~. \label{fDM2}
\end{eqnarray}
Note that the upright "e" represents the base of the natural logarithm, while the black body "$\bfe_z$" denotes the basis vector of the $z$ coordinate.

For all cases, all the parameters are normalized values with units defined in Sec.\ref{The FPR collision equation}. Unless otherwise specified, the default values for $(n_2, N_0)$ are set to $(7,7)$ as specified in section\SEC{Discretization of speed coordinate} and initial number of King functions $\NKo$ is set to 1 at the initial time level. Hence, default values of parameters are $\nhas=\vhaths=1$ and $\uzhas=\uzha$. For L01jd2 scheme, number of King function is constant, $\NKa \equiv \NKo$, while according to limitations provided by Algorithm \ref{alg: NKa}, $\NKmin \le \NKa(t) \le \NK \le \NKmax$, is specified in L01jd2NK and L01jd2nh schemes. The details are provided in Sec.\SEC{Update NK}. In this paper, the default solver is L01jd2NK scheme for scenarios where $max(m_M, 1/ m_M) \sim 1$, and L01jd2 scheme for cases where $max(m_M, 1/ m_M) \gg 1$ during collisions between electrons and ions.

\subsection{Two-species thermal equilibration}
\label{Two-species thermal equilibration}
 
In this instance, we demonstrate the convergence performance of our algorithm on two-species thermal equilibration, a widely used benchmark for evaluating schemes to solve the FPRS collision equation. The parameters for this case are $m_a=2$,  $m_b=3$, $Z_a=Z_b=1$, $n_a=n_b=1$, $\uza=\uzb=0$, $T_a=10$ and $T_b=20$. Theoretically values for the temperature and momentum at equilibrium state can be obtained using conservation equations for momentum and energy:
\begin{eqnarray}
    \Izsk  &=& \sum_a \left(m_a \nak \uzak \right) \ \equiv \  \sum_a \left(m_a \nako \uzako \right), \label{Izsk}
    \\
    \Ksk  &=& \frac{1}{2} \sum_a \left[3 \nak \Tak + m_a \nak (\uzak)^2\right] \ \equiv \ \frac{1}{2} \sum_a \left[3n_a \Tako + m_a \nako (\uzako)^2 \right] ~. \label{Ksk}
\end{eqnarray}
The initial total momentum and energy are $\Izsko=0$ and $\Ksko=45$ respectively. According to the conservation equations\EQ{Izsk}-\EQo{Ksk}, the finial average velocity and temperatures of the thermal equilibrium state should be $u_\infty=0$ and $T_\infty=15$ respectively.

The zeroth-order two-species thermal equilibration model, as described by Braginskii\cite{Braginskii1958TransportPlasma} , is presented through a semianalytical asymptotic temperature evolution equation:
\begin{eqnarray}
    \partial_t  T_a  &=&  \nuTab \left (T_a - T_b \right), \label{nufMab}
\end{eqnarray}
where the characteristic collision rate $\nuTab$ is given by\cite{Huba2011NRLFORMULARY} :
\begin{eqnarray}
    \nuTab  &\approx& 441.72 \times \frac{\sqrt{m_a m_b} \left(Z_a Z_b \right) ^2 n_b}{\left(m_a T_b + m_b T_a \right) ^{3/2}} \ln{\Lambda_{ab}} ~. \label{nuTab}
\end{eqnarray}
The temperature relaxation time is defined as $\tauTab = 1/ \nuTab$. The characteristic time $\tau_0$ is equivalent to the initial temperature relaxation time unless otherwise specified.
 
 \begin{figure}[H]
	\begin{center}
		\includegraphics[width=0.65\linewidth]{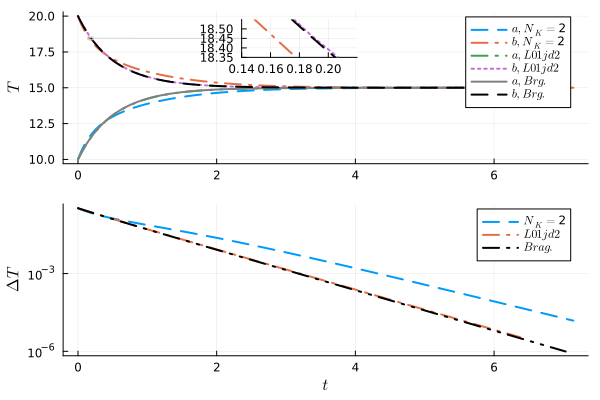}
	\end{center}
	\caption{Two-species temperature equilibration with fixed timestep, $\Delta t=2^{-5}$: Temperatures (upper) and relative disparity of temperatures (lower), $\Delta T  = \left |\Tak - \Tbk \right| / \left (\Tak + \Tbk \right) \label{DTabk}$ as functions of time. The color lines represent the numerical solution while the gray and black lines represent the Braginskii's values.}
	\label{FigTaTb2}
  \end{figure} 

The semianalytical equation\EQ{nufMab} is solved employing the standard explicit Runge-Kutta\cite{Rackauckas2017DifferentialEquations.jlaJulia} method of order 4, and the results are depicted in Fig.\FIG{FigTaTb2}. The temperatures are plotted as functions of time and compared against the numerical solution of our kinetic model with a fixed timestep of $\Delta t = 2^{-5}$ and L01jd2NK scheme with the maximum number of King function $\NK=2$. Fig.\FIG{FigTaTb2} demonstrates excellent agreement between our fully kinetic model and the semianalytical solution. Furthermore, upon comparing the results of L01jd2 ($\NKa(t) \equiv 1$ in this case) and L01jd2NK ($\NK=2$) with the semianalytical solution, it is observed that the temperature decay rate of L01jd2NK is a slightly faster than that of Braginskii's within the initial time block ($t < 0.58$ in this case), but this trend reverses when $t \ge 0.58$ (Similar to the behavior observed in the FVM's approach as shown in the Fig.14 of Ref.\cite{Taitano2015AEquation}). However, the results of L01jd2 scheme strictly adhere to the semianalytical solution, as shown in the lower subplot of Fig.\FIG{FigTaTb2}.
  
\begin{figure}[H]
	\begin{center}
		\includegraphics[width=0.65\linewidth]{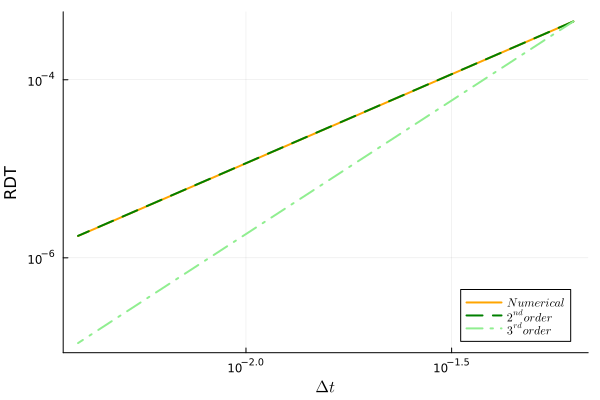}
	\end{center}
	\caption{Two-species temperature equilibration with fixed timestep, $\Delta t$: Demonstration of second-order convergence of the time discretization scheme.}
	\label{FigRDT2th}
  \end{figure}

The solid line in Fig.\FIG{FigRDT2th} depicts the relative deviation of temperature between our kinetic model and reference values $T_{ref}$, given by $RDT  = \left |T_{kin} - T_{ref} \right| / T_{ref} $ at $t=0.5$, as a function of the fixed timestep, $\Delta t$. The reference value, $T_{ref}$, is computed using a sufficiently small timestep, $\Delta t = 2^{-11} $, in our kinetic model. The results indicate that our algorithm exhibits $2^{nd}$-order convergence in time discretization.

\begin{figure}[H]
	\begin{center}
		\includegraphics[width=0.65\linewidth]{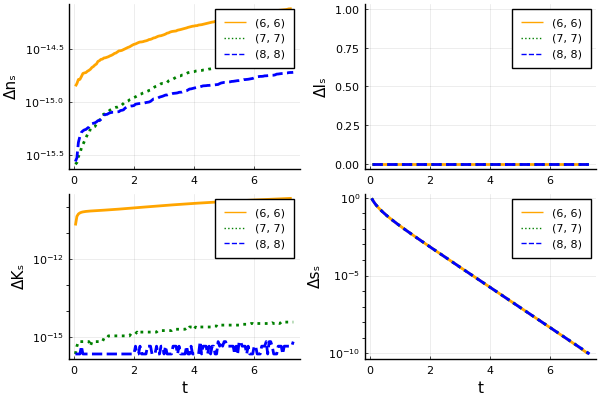}
	\end{center}
	\caption{Two-species temperature equilibration without enforcing conservation: Discrete conservation errors as functions of time and number of
    field nodes, $(n_2,N_0)$.}
	\label{FigCnIKsM}
  \end{figure}

The temporal evolution of errors in the conserved quantities, specifically discrete number (or mass) density, momentum, energy conservation,
\begin{eqnarray}
    \Delta \ns (\tk)  &=& (N_s)^{-1} \sum{_a} \left |\nak / \nako - 1 \right|, \label{Dns}
    \\
    \Delta \Izs (\tk)  &=&  \left |\Izsk / \Izsko - 1 \right|, \label{DIzs}
    \\
    \Delta \Ks (\tk)  &=&  \left |\Ksk / \Ksko - 1 \right|, \label{DKs}
\end{eqnarray}
and entropy conservation,
\begin{eqnarray}
    \Delta \ss (\tk)  &=& \left(\ssk - s_s^{k-1} \right) / \left [\Delta t   \left(\ssk + \delta s_s \right) \right] \label{Dtss}
\end{eqnarray}
are depicted in Fig.\FIG{FigCnIKsM} for varying number of field nodes with $(n_2,N_0)$. Here, the parameter $\delta s_s$ is given by $\delta s_s = \delta_{\zeta}^{-1} (\left |\ssko\right| + \left |\ss^{end}\right|)$, where $\zeta =  \left |\ssko \ss^{end}\right|/(\ssko \ss^{end})$ and $\ss^{end}$ are the value of $\ss$ at the finial time level. $\delta \ss \equiv 0$ in all cases except for the last one in Sec.\SEC{Numerical results}. Since entropy conservation serves as a convergence criterion for our algorithm, a first-order implicit scheme is utilized to calculate the entropy change, as defined by Eq.\EQ{Dtss}.

The discrete mass conservation is achieved with high precision for all given $(n_2,N_0)$, even without enforcing conservation. The discrete momentum conservation and H-theorem are preserved all the time. The discrete error of the energy conservation rapidly decreases with an increase of number of field nodes and reaches the level of round-off error when $(n_2,N_0)=(8,8)$, corresponding to a total number of nodes $N_n=257$. Since $\lMI=1$ for two-species temperature equilibration, the total number of field nodes is $N_{\rmv 2}=\lMI \times N_{\rmv}=1793$ for Rosenbluth potentials step \ref{Rosenbluth potentials}. The results of $\Delta \Ks$ in Fig.\FIG{FigCnIKsM} indicate that the convergence order of the velocity-space discretization scheme is about 16. The convergence order\cite{Fornberg1996AMethods} of a discrete algorithm is evaluated using:
\begin{eqnarray}
    order  &=& \left[1 / \log(2) \right] {\log(\epsilon_{n_2-1}/\epsilon_{n_2})}
    ~. \label{orderConverg}
\end{eqnarray}
Here, we assume that $\epsilon_{n_2} = \scrO [(h_{n_2})^{order}]$, where $h_{n_2}$ represents the grid size with grid number is $2^{n_2}+1$. The number of time levels in this case is determined to be $N_t=122$, based on the termination condition $\Delta s_s \le 10^{-10}$.

  \begin{figure}[H]
	\begin{center}
		\includegraphics[width=0.7\linewidth, height=0.29\linewidth]{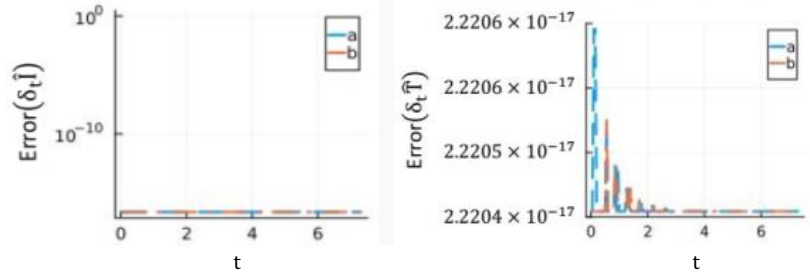}
	\end{center}
	\caption{Two-species temperature equilibration without enforcing conservation: Romberg integral errors of $\pDt \Ih$ and convergence criterion $\pDt \Th$ respect to time $t$ when $(n_2,N_0)=(7,7)$.}
	\label{FigCerror77}
  \end{figure}
  
\begin{figure}[H]
	\begin{center}
		\includegraphics[width=0.67\linewidth]{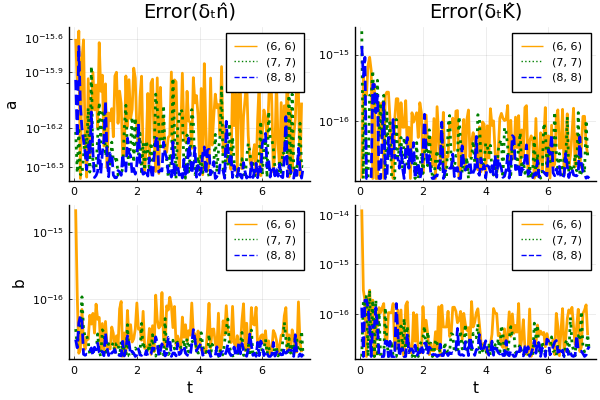}
	\end{center}
	\caption{Two-species temperature equilibration without enforcing conservation: Romberg integral errors of $\pDt \nh$ and $\pDt \Kh$ respect to time for collision species $a$ and $b$ with various $(n_2,N_0)$.}
	\label{FigCerrorTaTbM}
  \end{figure}

The relative errors of Romberg integrals for the first few orders of $\calRhjl$\EQ{Rhjl}, $\pDt \nha$\EQ{Rhna}, $\pDt \Izha$\EQ{RhIa}, $\pDt \Kha$\EQ{RhKa} and convergence criterion $\pDt \Tha$\EQ{RhTa} during two-species Coulomb collisions are depicted in Fig.\FIG{FigCerror77} and Fig.\FIG{FigCerrorTaTbM} when $(n_2,N_0)=(7,7)$. It is evident that all the relative errors are at the level of round-off error. The maximal errors occur at the initial moment and diminish to the level of round-off error over a collision time scale. The fact that convergence criterion $\pDt \Tha$ consistently equals the theoretical value indicates that our algorithm exhibits strong convergence properties.
  
The relationship between relative errors and orders $(n_2,N_0)$ is depicted in Fig.\FIG{FigCerrorTaTbM}. In all cases, the relative errors are significantly smaller that one,
indicating that condition\EQgive{Ccond} is consistently satisfied. Specifically, the relative errors of species $a$ remain at the level of round-off error throughout. However,  the orange solid line representing species $b$ in Fig.\FIG{FigCerrorTaTbM} reveals that both the relative errors of $\pDt \nhb$ and $\pDt \Khb$ are two orders of magnitude larger than the round-off error when $(n_2,N_0)=(6,6)$ at initial moments. This discrepancy can result in discrete errors in energy conservation (see Fig.\FIG{FigCnIKsM}). Furthermore, Fig.\FIG{FigCerrorTaTbM} also indicates that as the number of grids $(n_2,N_0)$ increases, the errors of species $b$ rapidly decrease to the level of round-off error. This suggests that integrals\EQ{Rhna}-\EQo{RhKa} will be accurate at least for one species during two-species collision process and the integration accuracy for any less accurate species can be improved by refining field nodes. 

  \begin{figure}[H]
	\begin{center}
		\includegraphics[width=0.67\linewidth]{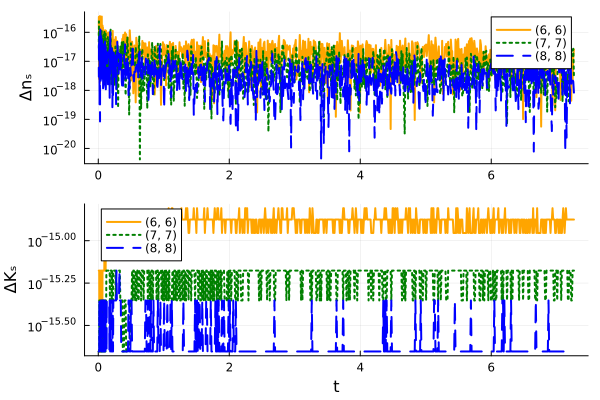}
	\end{center}
	\caption{Two-species temperature equilibration with enforcing conservation: Discrete conservation errors as functions of time and number of field nodes, $(n_2,N_0)$.}
	\label{FigCnKsMC2}
  \end{figure}

When implementing the conservation enforcing algorithm\SEC{alg: Conservation enforcing}, the time histories of the associated conservation errors, $\Delta n_s$ and $\Delta K_s$ are plotted in Fig.\FIG{FigCnKsMC2} for varying number of field nodes, $(n_2,N_0)$. As anticipated, precise discrete conservation can be achieved for all given $(n_2,N_0)$.
 
\begin{figure}[H]
	\begin{center}
		\includegraphics[width=0.67\linewidth]{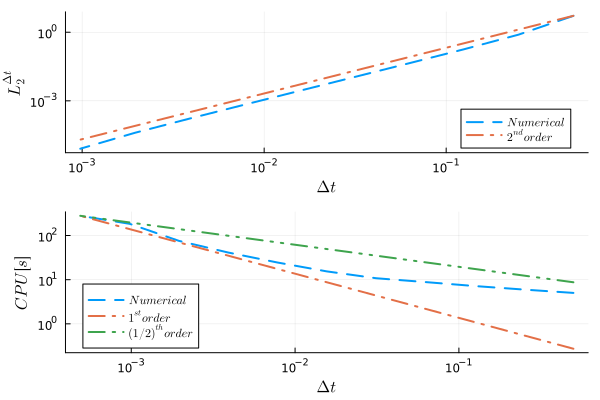}
	\end{center}
	\caption{Two-species temperature equilibration with fixed timestep, $\Delta t$: Demonstration of $2^{th}$-order convergence of time discretization scheme in $L_2$-norm and the CPU time as function of timestep, $\Delta t$.}
	\label{FigL2dt}
  \end{figure}

We can also demonstrate second-order time convergence of the trapezoidal scheme by computing the $L_2$-norm of relative difference in solution compared to a reference solution,
\begin{eqnarray}
    L_2^{\Dt}  &=& \sqrt{\left < \fl{^{\Delta t}} / \fl{^{\Delta t_{ref}}} - 1, \fl{^{\Delta t}} / \fl{^{\Delta t_{ref}}} - 1 \right >} ~. \label{L2Dt}
\end{eqnarray}
Here, $\fl{^{\Delta t_{ref}}}$ is the solution obtained using a reference timestep size, $\Delta t_{ref} = 2^{-11}$; refer to Fig.\FIG{FigL2dt} (upper) when $(n_2,N_0)=(7,7)$. As expected, second-order convergence is realized with the refinement of $\Dt$. The CPU time as a function of $\Delta t$ when $(n_2,N_0)=(7,7)$ is also plotted in the lower subplot of Fig.\FIG{FigL2dt}. The total solution time scales approximately as $\scrO (1 / \Delta t)$. This first-order correlation arises from the rapid convergence of implicit iteration, typically requiring  only a few iterations even with different timesteps , as demonstrated by Fig.\FIG{FigtNtRDt} in Sec.\SEC{Three-species e-D-A thermal equilibration}. Compared to the explicit timestep, $\Delta _t^{Exp} \approx 1.5 \times 10^{-3}$  in this case estimated by Eq.\EQ{dtCFL}, a timestep greater than one or two order of magnitude can be used in our algorithm with acceptable time-discrete precision. 
  
\begin{figure}[H]
	\begin{center}
		\includegraphics[width=0.7\linewidth]{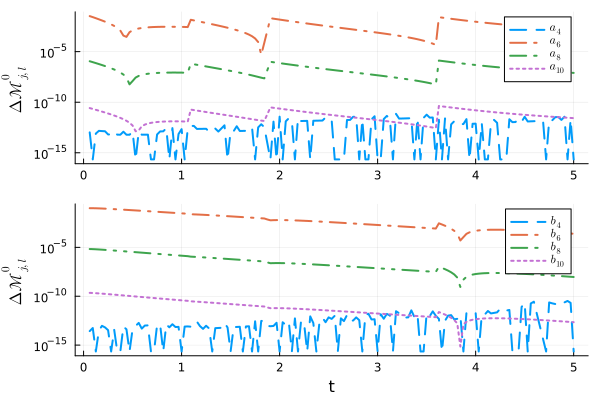}
	\end{center}
	\caption{Two-species temperature equilibration without enforcing conservation: Time discretization errors of the non-conserved moments as functions of time with various $j$ when $\Delta t = 2^{-5}$, $\NK=2$ and $(n_2,N_0)=(7,7)$.}
	\label{FigerrMhcopM4}
  \end{figure}

Fig.\FIG{FigerrMhcopM4} illustrates the time discretization errors of the non-conserved moments, $\Delta \calMhjl$ as functions of time with various $j$ when $j \ge 4$. Here, $\Delta t = 2^{-5}$, $\NK=2$ and $(n_2,N_0)=(7,7)$. The quantity $\Delta \calMhjl$, which is more effective than $\delta \calMhjl$ presented by Eq.\EQ{dMhjl}, is defined as:
\begin{eqnarray}
    \Delta \calMhjl  &=& \left |\hat{M}_{j,l}(\tkI) - \calMhjl(\tkI) \right| / \left |\calMhjl(\tkI) - \hat{M}_{j,l}(\tk) \right|, \label{DMhjl}
\end{eqnarray}
which measures the relative error resulting from velocity discretization during optimization process. Notes that the symbols $\hat{M}_{j,l}$ and $\calMhjl$ represent the normalized kinetic moment, which are computed from the amplitude function before\EQgive{Mhjl} and after\EQgive{CPEs} being smoothed by King function, respectively. By applying Eqs.\EQ{uzhanuTs}-\EQo{KhanuTs} in L01jd2NK scheme, the time discretization errors of the conserved moments, $\Delta \calMhjl$ given by Eq.\EQ{DMhjl}, for all species are exactly zero. Furthermore, convergence of the optimization of $(j,l)^{th}$-order normalized kinetic moments is achieved when $\Delta \calMhjl \le rtol_{NK}$, where $rtol_{NK}$ is a given relative tolerance. In this study, we set parameter $rtol_{NK}$ to $10^{-11}$ unless otherwise specified. As shown in Fig.\FIG{FigerrMhcopM4}, the moment with order $j=6$ exhibits a maximum deviation not exceeding $6.6\%$ for all species under consideration when $\NK=2$. Additionally, the time discretization errors for the convergent order, $j=4$, are generally no greater than $10^{-11}$.

  \begin{figure}[H]
	\begin{center}
		\includegraphics[width=0.7\linewidth]{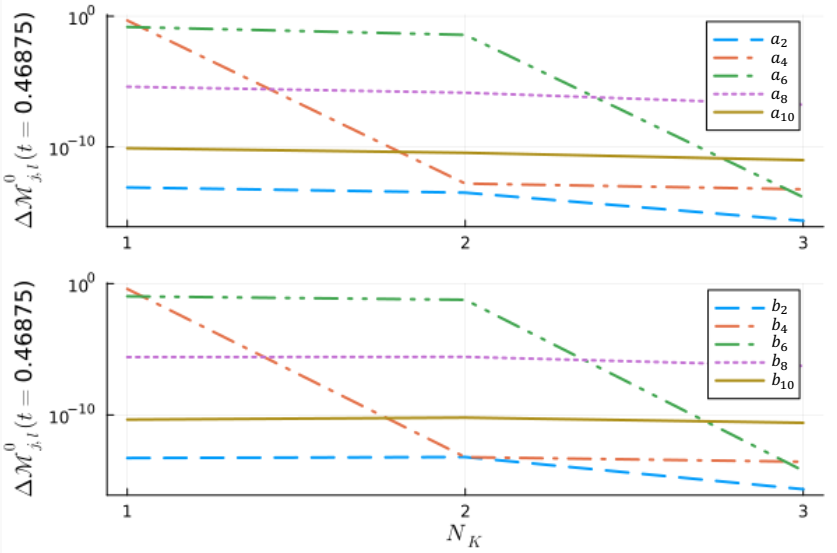}
	\end{center}
	\caption{Two-species temperature equilibration without enforcing conservation: Time discretization errors of the non-conserved moments as functions of $\NK$ with various $j$ when $\Delta t = 2^{-6}$, $\tk\approx0.47$ and $(n_2,N_0)=(8,8)$.}
	\label{FigerrMhcopM4NK}
  \end{figure}

The high-order moment convergence property of the present method is further investigated. Fig.\FIG{FigerrMhcopM4NK} illustrates the time discretization errors of the non-conserved moments, $\Delta \calMhjl$, at $\tk\approx0.47$ as functions of $\NK$ with various $j$. For this test, a refined timestep and field nodes, $\Delta t = 2^{-6}$ and $(n_2,N_0)=(8,8)$ are utilized. As depicted in Fig.\FIG{FigerrMhcopM4NK}, all the discretization errors of the non-conserved moments, $\Delta \calMhjl$, decrease as $\NK$ increase. Furthermore, it is observed that the highest convergent order is $j=2$ when $\NK=1$, $j=4$ when $\NK=2$ and $j=6$ when $\NK=3$. Therefore, it can be concluded that King method (Sec.\SEC{King method}) is a moment convergence algorithm. 

  \begin{figure}[H]
	\begin{center}
		\includegraphics[width=0.65\linewidth]{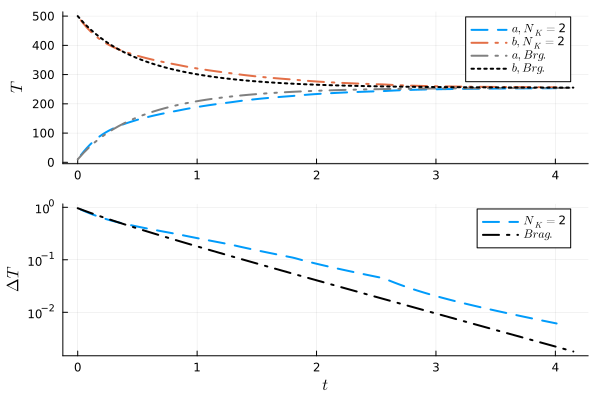}
	\end{center}
	\caption{Two-species temperature equilibration with temperatures, $T_a=10$, $T_b=500$ and fixed timestep, $\Delta t=2^{-5}$: When $\NK=2$ and $(n_2,N_0)=(8,7)$.}
	\label{FigTaTb2Ta10Tb500}
  \end{figure}
  
\begin{figure}[H]
	\begin{center}
		\includegraphics[width=0.65\linewidth]{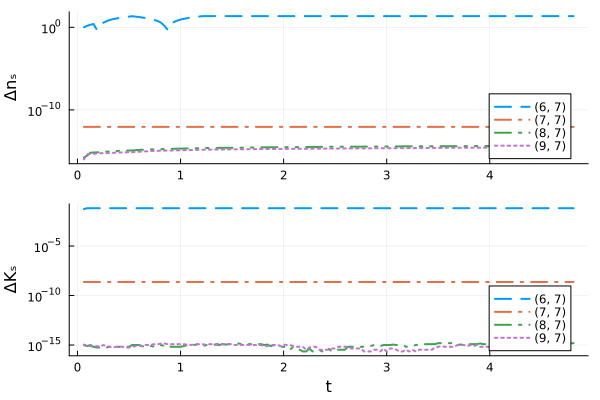}
	\end{center}
	\caption{Two-species temperature equilibration with temperatures, $T_a=10$, $T_b=500$ and fixed timestep, $\Delta t=2^{-5}$: Discrete conservation errors as functions of time and $n_2$ when $N_0=7$.}
	\label{FigCnKMTa10Tb500}
  \end{figure}

In situations with large thermal velocity disparity where $v_{th_f}/v_{th_s} \gg 1$, the maximum value of normalized speed $\rmvh_{M}$ at $(k+1)^{th}$ time level within $(k+1)^{th}$ time block is determined by a more sophisticated method, for species $a$, can be described as follow:
\begin{eqnarray}
    \delta_l^0 \left(\rmvh_M^{k+1} \right)^{j} \fhl (\rmvh_M^{k+1}, \tkI) + \delta_l^1 \left(\rmvh_M^{k+1} \right)^{j+1} \left|\fhl (\rmvh_M^{k+1}, \tkI) \right|  &=& 4.44 \times 10^{-17}, \ j = 2.5 \NKa. \quad \label{vhmaxAdapt}
\end{eqnarray}
Applying information at the $(k+1)^{th}$ time level within $k^{th}$ time block and solving the above equation, the optimized value of  $\rmvh_{M}$ at $(k+1)^{th}$ time level within $(k+1)^{th}$ time block can be obtained. It is evident that $\rmvh_{M}$ will vary across different time levels, while remaining constant within each time block. This step ensures that the value of distribution function at the right boundary of speed can be effectively disregarded as zero.

For example, when $T_a=10$ and $T_b=500$, with all other parameters remaining the same as in the previous case, the temperatures calculated by the L01jd2NK scheme using a fixed timestep of $\Delta t=2^{-5}$, $\NK=2$ and $(n_2,N_0)=(8,7)$ are depicted as functions of time in Fig.\FIG{FigTaTb2Ta10Tb500}. As anticipated, the results show good agreement between our fully kinetic model and the semianalytical solution. By comparing Fig.\FIG{FigCnKMTa10Tb500} and Fig.\FIG{FigCnIKsM}, we can see that achieving the same level of precision with an increased temperature difference requires more refined field nodes by increasing $n_2$.

\subsection{Electron-Deuterium thermal equilibration}
\label{Electron-Deuterium thermal equilibration}

In order to verify the convergence in scenarios involving significant mass disparity, we examine the thermal equilibration between electron and Deuterium, denoted as slanted $e$ and $D$ respectively. The parameters $m_e=1/1836$, $m_D=2$, $-Z_e=Z_D=1$, $n_e=n_D=1$, $\uh_e=\uh_D=0$, $T_e=1$ and $T_D=10$. In this instance, the final average velocity and temperatures of the thermal equilibrium state are expected to be $u_\infty=0$ and $T_\infty=5.5$ respectively. The L01jd2 scheme is employed for solving this case with a self-adaptive timestep. The total number of time levels, denoted as $N_t$, is $792$,  when the termination condition is $\Delta s_s \le 10^{-8}$. We aim to demonstrate that L01jd2 scheme provides a reliable approximation of the fully kinetic model for situations characterized by large mass disparity.

  \begin{figure}[H]
	\begin{center}
		\includegraphics[width=0.65\linewidth]{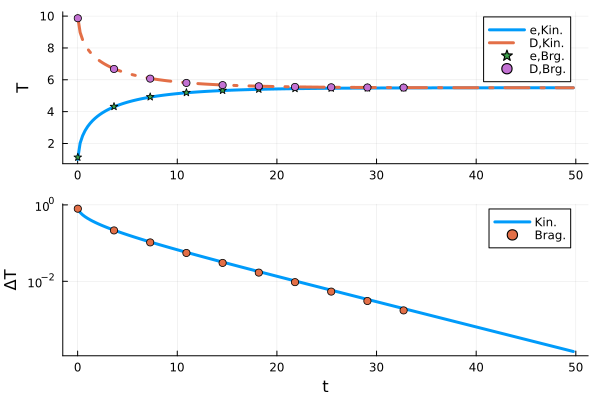}
	\end{center}
	\caption{$e$-$D$ thermal equilibration with enforcing conservation when solved by L01jd2 with fixed timestep, $\Delta t = 2^{-5}$: Average velocities and temperatures as function of time $t$ when $(n_2,N_0)=(7,7)$.}
	\label{FigDTeD}
  \end{figure}

\begin{figure}[H]
	\begin{center}
		\includegraphics[width=0.63\linewidth]{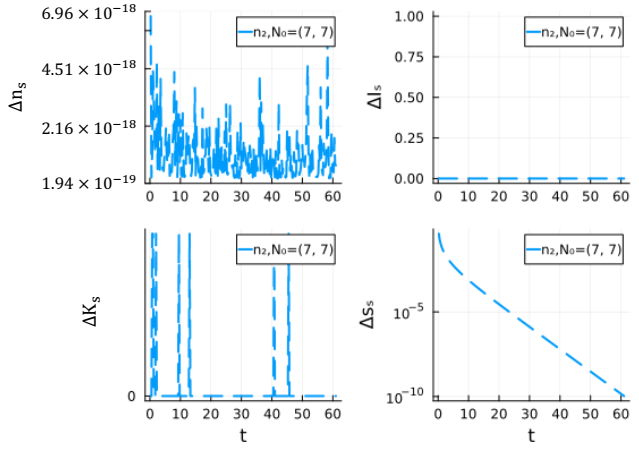}
	\end{center}
	\caption{$e$-$D$ thermal equilibration with enforcing conservation when solved by L01jd2 with fixed timestep, $\Delta t = 2^{-5}$: Discrete conservation errors as functions of time when $(n_2,N_0)=(7,7)$.}
	\label{FigeDCnIKsfM}
  \end{figure}

\begin{figure}[htp]
	\begin{center}
		\includegraphics[width=0.63\linewidth]{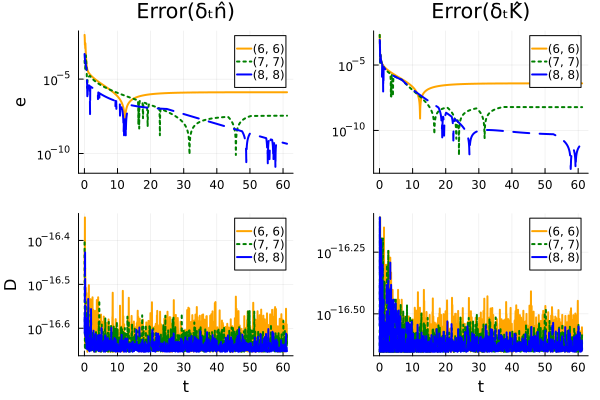}
	\end{center}
	\caption{$e$-$D$ thermal equilibration with enforcing conservation when solved by L01jd2 with fixed timestep, $\Delta t = 2^{-5}$:  Local errors of $\pDt \nh$ and $\pDt \Kh$ as functions of time with various $(n_2,N_0)$.}
	\label{FigeDCerrfM}
  \end{figure}

The temperatures of species $e$ and $D$ are presented as functions of time $t$ in Fig.\FIG{FigDTeD}, demonstrating the attainment of correct equilibrium values. The temporal evolution of errors in discrete number density, momentum, energy conservation\EQ{Dns}-\EQo{DKs} and entropy conservation is illustrated in Fig.\FIG{FigeDCnIKsfM}. Additionally, the local relative errors, $Error(\pDt \nh)$ and $Error(\pDt \Kh)$ of species $e$ and $D$ are depicted in Fig.\FIG{FigeDCerrfM} when enforcing conservation. It can be observed that all errors of species $D$ are at the level of round-off errors, which those for species $e$ are acceptably small, aligning with the convergence criterion for conservation\EQgive{Ccond}.

To assess the effectiveness of L01jd2 in scenarios involving significant mass disparity, we validate its convergence with the following criterion:
\begin{eqnarray}
     \Delta_0 \fl(\rmvhA,\tkI) &=& \left|\fl{^*}(\rmvhA,\tkI) - \fl(\rmvhA,\tkI)  \right| / \left| \fl{^*}(\rmvhA,\tkI) - \fl(\rmvhA,\tk) \right| = C_{0}  \label{C0}
\end{eqnarray}
and
\begin{eqnarray}
     \Delta_2 \fl &=& \left({N_v}\right)^{-1} \sum{_{\alpha=1}^{N_v} }\delta \fl(\rmvhA,\tkI)  ~. \label{C2} 
\end{eqnarray}
Here,
\begin{eqnarray}
     \delta \fl(\rmvhA,\tkI) &=& \left|\fl{^*}(\rmvhA,\tkI) - \fl(\rmvhA,\tkI)  \right| / \left|\Delta t \times \fl{^*}(\rmvhA,\tkI) \right| = C_2 (\Delta t)^2 ~.
\end{eqnarray}
Functions $\fl{^*}(\rmvh,\tkI)$ and $\fl(\rmvh,\tkI)$ denote amplitudes before and after being smoothed by King function at the $(k+1)^{th}$ time level within $k^{th}$ time block for species $a$, respectively.

\begin{figure}[htp]
	\begin{center}
		\includegraphics[width=0.65\linewidth]{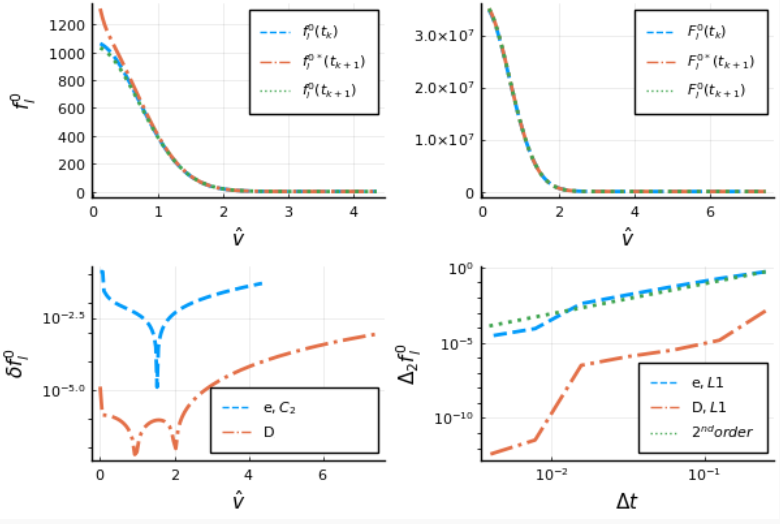}
	\end{center}
	\caption{e-D temperature equilibration with enforcing conservation when solved by L01jd2 with fixed timestep, $\Delta t = 2^{-5}$: Criterion $C_2$\EQ{C2} for situation with large mass disparity.}
	\label{FigDTRDKingC2}
  \end{figure}

The convergence of function $\Delta_2 \fl$ is illustrated in Fig.\FIG{FigDTRDKingC2}, with a fixed finial time of $\tkI=0.5$. The right-lower subplot in Fig.\FIG{FigDTRDKingC2} presents compelling evidence of the second-order accuracy convergence convergence of King method as $\Delta t$ is refined. When using a sufficiently small timestep, such as $\Delta t = 0.01$, the maximum relative disparity of the distribution function before and after being smoothed by the King function does not exceed $11\%$ for species $e$ and $10^{-5}$ for species $D$ in this scenario. This can be observed from the distribution function at $k^{th}$ and $(k+1)^{th}$ time levels for species $e$ and $D$ in Fig.\FIG{FigDTRDKingC2} (upper). The detailed relative disparity as a function of $\rmvh$ is plotted in the left-lower subplot of Fig.\FIG{FigDTRDKingC2} when $\Delta t = 2^{-5}$.

\begin{figure}[H]
	\begin{center}
		\includegraphics[width=0.65\linewidth]{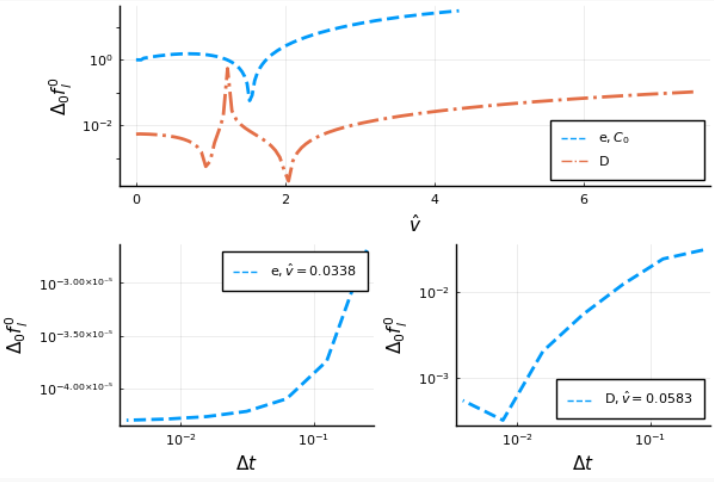}
	\end{center}
	\caption{e-D temperature equilibration with enforcing conservation when solved by L01jd2 with fixed timestep, $\Delta t = 2^{-5}$: Criterion $C_0$\EQ{C0} for situation with large mass disparity.}
	\label{FigDTRDKingC0}
  \end{figure}

Similarly, the convergence of function $\Delta_0 \fl$ at the grid point $v_{\alpha}=1.68\times 10^{-2}$ is demonstrated in Fig.\FIG{FigDTRDKingC0}. The lower subplots show that the parameter $C_0$ tends to  stabilize as a constant with the refinement of $\Delta t$ at single speed node, i.e. $\rmvhA=0.0338$ for species $e$ and $\rmvhA=0.0583$ for species $D$. For species $D$, which exhibits a higher precision of distribution function, $\Delta_0 \fl(\rmvhA=0.0583)$ remains bellows $0.1$ for all provided timesteps and does not exceed $1\%$ when $\Delta t = 2^{-5}$. Considering the high precision of conservation in discrete with a acceptable timestep $\Delta t = 2^{-5}$, the L01jd2 scheme is a good approximation of the fully kinetic model for situations with large mass disparity.
  
\subsection{Electron-Deuterium temperature and momentum equilibration}
\label{Electron-Deuterium disparate temperature and momentum equilibration}

We consider $e$-$D$ temperature and momentum equilibration with parameters $m_e=1/1836$, $m_D=2$, $-Z_e=Z_D=1$, $n_e=n_D=1$, $\uh_e=0.1$, $\uh_D=-3.162\times 10^{-2}$, $T_e=1$ and $T_D=100$ when $(n_2,N_0)=(7,7)$. The initial value of $\lMI=13$, $\veth\approx0.063$ and $\vDth\approx0.010$. The expected finial average velocity and temperatures are approximately $u_\infty \approx -9.769\times 10^{-2}$ and $T_\infty \approx 50.504$, respectively. This case will also be solved by the L01jd2 scheme with a self-adaptive timestep, where the number of time levels is determined to be $N_t=4877$. In theory, electron-deuterium momentum equilibration will occur first followed by reaching a state of temperature equilibration.

  \begin{figure}[H]
	\begin{center}
		\includegraphics[width=0.6\linewidth]{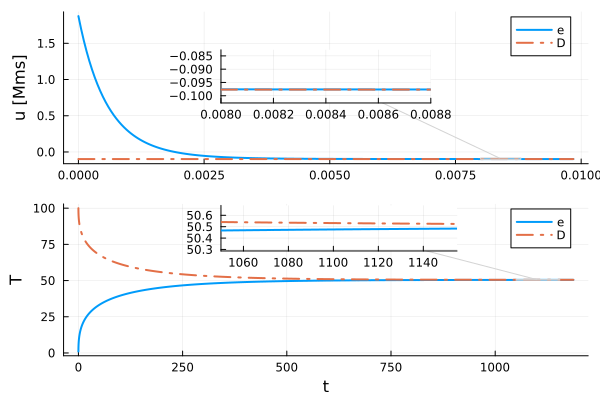}
	\end{center}
	\caption{$e$-$D$ thermal and momentum equilibration with enforcing conservation: Average velocities and temperatures as function of time $t$ when $(n_2,N_0)=(7,7)$.}
	\label{FigDuTeD}
  \end{figure}
  
The average velocities and temperatures of species $e$ and $D$ are depicted as functions of time $t$ in Fig.\FIG{FigDuTeD}. The lens in Fig.\FIG{FigDuTeD} demonstrates that correct equilibrium values of momentum and temperature are achieved. As anticipated, the characteristic time for momentum relaxation time is significantly shorter than the characteristic time for temperature relaxation during the process of electron-deuterium temperature and momentum equilibration. 
  
  \begin{figure}[H]
	\begin{center}
		\includegraphics[width=0.63\linewidth]{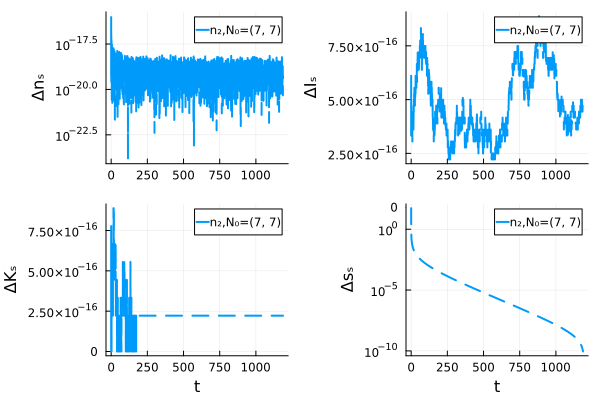}
	\end{center}
	\caption{$e$-$D$ thermal and momentum equilibration with enforcing conservation when $T_e=1$ and $T_D=100$: Discrete conservation errors as functions of time when $(n_2,N_0)=(7,7)$.}
	\label{FigeDCnIKTa1Tb100}
  \end{figure}

  \begin{figure}[H]
	\begin{center}
		\includegraphics[width=0.63\linewidth]{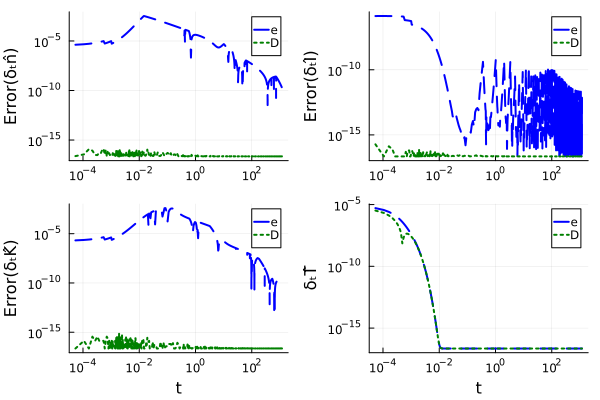}
	\end{center}
	\caption{$e$-$D$ thermal and momentum equilibration with enforcing conservation when $T_e=1$ and $T_D=100$:  Local errors of $\pDt \nh$, $\pDt \Ih$, $\pDt \Kh$ and $\pDt \Th$ as functions of time when $(n_2,N_0)=(7,7)$.}
	\label{FigeDCerrTa1Tb100M}
  \end{figure}

The time histories of the errors in discrete number density, momentum, energy conservation and entropy conservation are depicted in Fig.\FIG{FigeDCnIKTa1Tb100}. As before, mass, momentum and energy conservation\EQ{Dns}-\EQo{DKs} are enforced to the level of round-off error and H-theorem are preserved all the time, as demonstrated in Fig.\FIG{FigeDCnIKTa1Tb100}. Fig.\FIG{FigeDCerrTa1Tb100M} illustrates that the local relative errors, $Error(\pDt \nh)$, $Error(\pDt \Ih)$, and $Error(\pDt \Kh)$ are sufficiently small, satisfying the convergence criterion for conservation\EQ{Ccond} during $e$-$D$ collision when $(n_2,N_0)=(7,7)$. As expected, the local relative errors of species $D$, with its larger mass, are smaller than those of species $e$. The convergence criterion $\pDt \Th$ (displayed in Eq.\EQ{RhTa}) is approximately valid with high precision, especially when $t \ge 10^{-2}$.

\subsection{Three-species (e-D-alpha) thermal equilibration}
\label{Three-species e-D-A thermal equilibration}

The finial test case involves the three-species ($e$-$D$-$alpha$) thermal equilibration, which is an crucial issue in fusion plasma. Theoretically, the high-energy $alpha$ particles in burning plasma will initially exchange energy with electrons (of comparable $\vth$) and later thermalize with $D$ particles. In this subsection, the subscript "$\alpha$" represents species $alpha$.

The simulation parameters are $m_e=1/1836$, $m_D=2$, $m_{\alpha}=4$, $-Z_e=Z_D=1$, $Z_{\alpha}=2$, $n_e=3$, $n_D=1$, $n_{\alpha}=1$. Initially $\uh_e=\uh_D=\uh_{\alpha}=0$, $T_e=T_D=1$ and $T_{\alpha}=1750$, we choose $(n_2,N_0)=(7,7)$. In this case, $\veth \approx 6.256 \times 10^{-2}$, $\vDth \approx 1.036 \times 10^{-3}$ and $\vAth \approx 3.065 \times 10^{-2}$. The characteristic time $\tau_0$ is equivalent to the initial temperature relaxation time between $e$ and $D$, $\tau_0= 1 / \nu_T^{eD}$, where 
\begin{eqnarray}
    \nu_T^{eD}  &\approx& 441.72 \times \frac{\sqrt{m_e m_D} \left(Z_e Z_D \right) ^2 n_D}{\left(m_e T_D + m_D T_e \right) ^{3/2}} \ln{\Lambda_{eD}} ~. \label{nuTeD}
\end{eqnarray}
The maximum thermal velocity ratio, $\veth / \vDth \sim 60$. 

  \begin{figure}[H] 
	\begin{center}
		\includegraphics[width=0.6\linewidth]{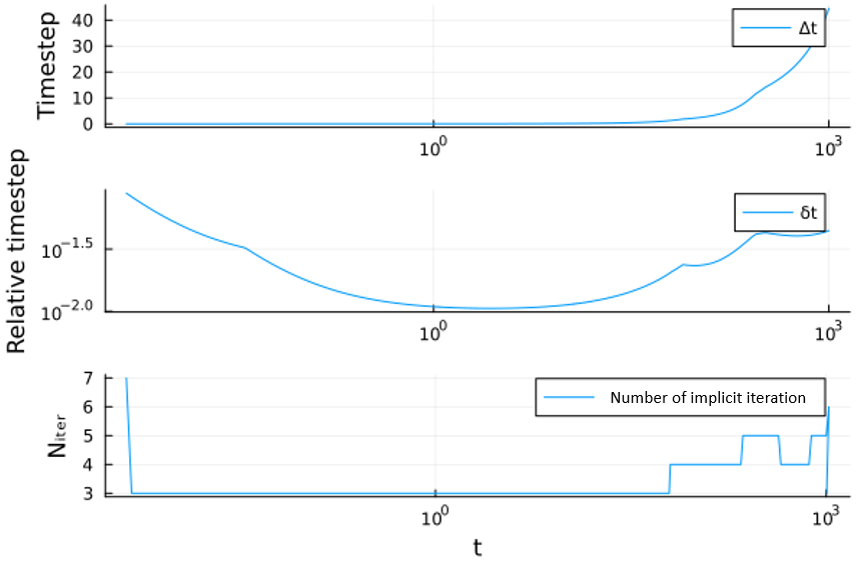}
	\end{center}
	\caption{$e$-$D$-$\alpha$ thermal equilibration: timestep $\Delta t$, the relative timestep $\delta t$ and number of implicit iterations, $N_{iter}$, within each time block as functions of time.}
	\label{FigtNtRDt}
  \end{figure}
  
This case is solved using a self-adaptive timesteps, with a total of 710 time levels. Fig.\FIG{FigtNtRDt} illustrates the timestep $\Delta t$, the relative timestep $\delta t$ and number of implicit iterations $N_{iter}$ within each time block as functions of time. Here, the $k^{th}$ timestep and relative timestep are respectively defined as $\Delta t_k = \tkI - \tk$ and $\delta t_k = \Delta \tk / \tkI$. The number of implicit iterations consistently remains below 10, regardless of the timestep and relative timestep. Additionally, due to the effectiveness of the King method, the value of $N_{iter}$ typically does not exceed 5 in most time block.

\begin{figure}[H] 
	\begin{center}
		\includegraphics[width=0.63\linewidth]{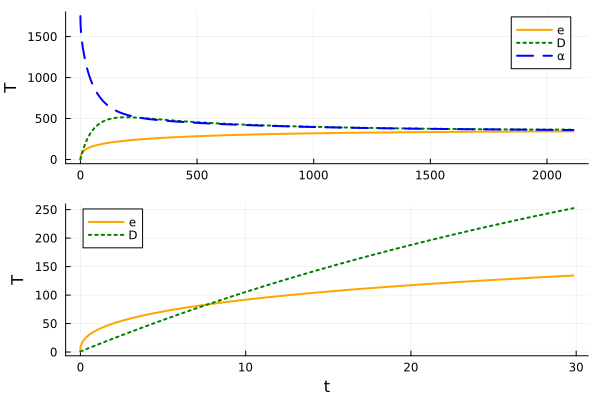}
	\end{center}
	\caption{$e$-$D$-$\alpha$ thermal equilibration: Temperature of the three species as functions of time (upper) and of $e$ and $D$ at the early time (lower).}
	\label{FigTeTDTA}
  \end{figure}

The temperatures of all species as functions of time $t$ are plotted in Fig.\FIG{FigTeTDTA}. As anticipated, the electrons exhibit a more rapid heating rate comparing to the $D$ particles during the early stage when $t \le 8$ (lower figure in Fig.\FIG{FigTeTDTA}). However, as electrons and $D$ quickly heat up while $\alpha$ particles cool down, the preferential interaction switches to one between $\alpha$ and $D$, ultimately leading to their thermalization together. Eventually, all three species reach the expected equilibrium temperature of $T_{\infty}=350.08$ (upper panel in Fig.\FIG{FigTeTDTA}).

  \begin{figure}[H]
	\begin{center}
		\includegraphics[width=0.63\linewidth]{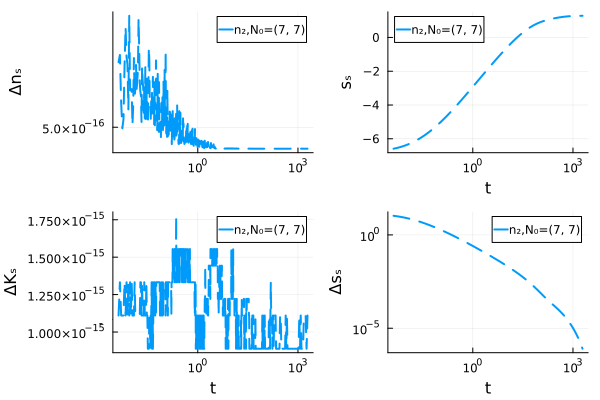} 
	\end{center}
	\caption{$e$-$D$-$\alpha$ thermal equilibration: Discrete conservation errors as functions of time.}
	\label{FigeDCnIKTe1TD1TA1750}
  \end{figure}

\begin{figure}[htp]
	\begin{center}
		\includegraphics[width=0.6\linewidth]{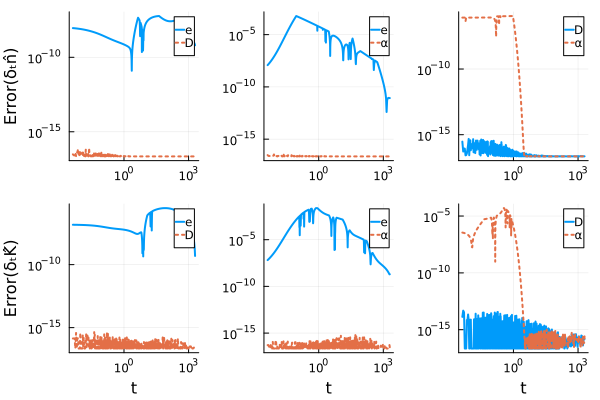}
	\end{center}
	\caption{$e$-$D$-$\alpha$ thermal equilibration: Local relative errors of $\pDt \nh$ and $\pDt \Kh$ during all the two-species collision processes as functions of time.}
	\label{FigeDACerrTeTDTA1750nK}
  \end{figure}

The temporal histories of local errors in discrete number density, momentum, energy conservation and entropy conservation are plotted in Fig.\FIG{FigeDCnIKTe1TD1TA1750}.
As expected, mass and energy conservation represented by Eq.\EQ{Dns} and Eq.\EQ{DKs} are enforced to the level of round-off error. The right-upper subplot of Fig.\FIG{FigeDCnIKTe1TD1TA1750} demonstrates that the entropy satisfies the H-theorem at all time. $\Delta s_s$\EQ{Dtss} is consistently non-negative, as showed in the right-lower subplot of Fig.\FIG{FigeDCnIKTe1TD1TA1750}. In this case, the value of $\delta \ss$ is about $-7.87$. The total local relative errors of $\pDt \nh$ and $\pDt \Kh$ for all species as functions of time $t$ are presented in Fig.\FIG{FigeDACerrTeTDTA1750nK}. It can be observed that the convergence criterion for conservation denoted by Eq.\EQ{Ccond} is also satisfied for all sub-processes involving two-species collision process.

\section{Conclusion}
\label{Conclusion}

In this study, a nonlinear framework (SHE together with KFE) has been introduced for solving the multi-species nonlinear 0D-2V axisymmetric Fokker-Planck-Rosenbluth (FPR) collision equation. In this framework, Legendre polynomial expansion is employed in the angular coordinate, which converges exponentially. KFE, a moment convergence method, is utilized in  speed coordinate. This approach ensures mass, momentum, energy conservation and satisfies the H-theorem for plasma simulations with general mass and temperature. An efficient implicit algorithm also has been constructed for weakly anisotropic plasmas based on this framework, employing the nonlinear Shkarofsky's formula of FPR (FPRS) collision operator. The time block technique (TBT) and moment optimization method is utilized in this nonlinear algorithm. Kinetic moments are computed by Romberg integration with high precision. Subsequently, post-step projection to manifold method is applied to enforce the exact conservation of the collision operators. 

The high accuracy of our algorithm is demonstrated by solving several typical problems in various non-equilibrium configurations. To handle the large disparate of thermal velocities resulting from the arbitrary disparity of mass and temperature, we incorporate  mapping between the different field nodes for collision species. This is also accomplished based on the King function. The fast convergence and high efficiency in handing various challenging problems have showcased the potential of our approach for multi-scale simulation of plasma.  

In order to fully realize the potential of the proposed framework for nonlinear, multi-scale plasma systems, it is necessary to expand our approach into a general self-adaptive scheme, including self-adaptive collection of $(j,l)$, number of King function $\NKa$ and field nodes determined by parameter $(n_2,N_0)$ in the speed coordinate. This will be discussed in a future publication. We finally remark that the SHE together with KFE framework is efficient for moderately anisotropic plasma systems, including weakly anisotropic plasmas, subsonic regions and low supersonic plasmas, while the L01jd2NK scheme is only designed for weakly anisotropic plasmas. The limitation for SHE together with KFE framework is the decreasing convergence rate of SHE as the ratio of average velocity to thermal velocity increases. While the limitation for L01jd2NK scheme is that characteristic parameter in the King function\EQ{King} may depend on $l$ when the characteristic group velocity, $\uzhas$, is large enough and tends to be 1. A general scheme based on the present framework for moderately anisotropic plasma systems will be addressed in our further research.

\section{Acknowledgments}
\label{Acknowledgments}

We express our gratitude to Peifeng Fan and Bojing Zhu for their valuable discussions. We particularly appreciate the assistance of Shichao Wu and Ye Zhu in enhancing the English language of our manuscript. This work is supported by the GuangHe Foundation (ghfund202202018672), Collaborative Innovation Program of Hefei Science Center, CAS, (2021HSC-CIP019), National Magnetic Confinement Fusion Program of China (2019YFE03 060000), Director Funding of Hefei Institutes of Physical Science from Chinese Academy of Sciences (Grant Nos. E25D0GZ5), and Geo-Algorithmic Plasma Simulator (GAPS) Project.

\appendix    

\begin{appendices}

\section{King Function}
\label{King Function}

When the velocity space exhibits axisymmetry with $\uha=\uzha \bfe_z$, the Gaussian function will be:
\begin{eqnarray}
    \scrG\left(\rmvh,\mu \right) &=& \frac{\nha}{(\vhath)^3} \e^{-\left [\left(\vh - \uzha \bfe_z \right )/\vhath  \right]^2}.  \label{GaussianFun}
\end{eqnarray}
The Gaussian function can be expanded using Legendre polynomials as:
  \begin{eqnarray}
      \scrG \left(\rmvh,\mu \right) &=& \sum{_{l=0}^{l_M}} \scrG_l \left( \rmvh \right) \Pl \left(\mu \right) ~. \label{GGl}
  \end{eqnarray}
The $l^{th}$-order amplitude, $\scrG_l \left( \rmvh \right)$ can be calculated by the inverse transformation of Eq.\EQ{GGl}, reads:
   \begin{eqnarray}
       \scrG_l \left(\rmvh \right) &=& \int_{-1}^1 \scrG\left(\rmvh,\mu \right)  \Pl (\mu) \rmd \mu, \quad l \in 0,1,\cdots,l_M ~. \label{GlGanalys}
   \end{eqnarray}
Substitute Eq.\EQ{GaussianFun} into the equation above, and after a tedious derivation process, one can obtain:
   \begin{eqnarray}
       \scrG_l \left(\rmvh \right) &=& \frac{\nha}{(\vhath)^3} \sum_{m=0}^{l} C_{G_l} \frac{1}{\xih_a^{m+1}} \left[(-1)^m \e^{-\left(\frac{\rmvh - \uzha}{\vhath} \right)^2} - (-1)^l \e^{-\left(\frac{\rmvh + \uzha}{\vhath} \right)^2} \right] ~. \label{GlG}
   \end{eqnarray}
Here, $\xih_a = 2 \rmvh \uzha / \vhath^2$ and
   \begin{eqnarray}
       C_{G_l} &=&  \frac{2l+1}{2} \frac{1}{2^m m!} \frac{(l+m)!}{(l-m)!} ~. \label{CGl}
   \end{eqnarray}
The $l^{th}$-order King function\EQ{King} is in direct proportion to $\scrG_l$, reads:
\begin{eqnarray}
    \Kl \left(\rmvh;\uzha,\vhath \right) &=& \sqrt{\frac{1}{2\pi}} \frac{1}{\nha} \scrG_l \left(\rmvh \right)~.  \label{KingGl}
\end{eqnarray}

The new function introduced in Sec.\SEC{King method}, King function\EQ{King} has the following properties. When $\rmvh \to \infty$, the King function will be:
\begin{eqnarray}
    \lim_{\rmvh \to \infty} \Kl \left(\rmvh;\iota,\sigma \right) &\to& 0 ~.  \label{Kinginfty}
\end{eqnarray}
When $\xi=2\iota \rmvh / \sigma^2 \to 0$, the King function has the following asymptotic behaviour:
\begin{eqnarray}
    \lim_{\xi \to 0} \Kl \left(\rmvh;\iota,\sigma \right) &=& \sqrt{\frac{1}{2\pi}} \frac{\e^{-\sigma^{-2} \left (\rmvh^2 + \iota^2  \right)}}{(2 l - 1)!!} \frac{\xi^l}{\sigma^3} \left [1 + \sum_{k=1}^{\infty} \frac{1}{2^l l!} \frac{(2 l - 1)!!}{(2 l + 2 k + 1)!!}  \xi^{2k} \right] ~.  \label{Kingxi}
\end{eqnarray}
In particular, $\xi \equiv 0$ gives:
\begin{eqnarray}
    \Kl \left(\rmvh;\iota,\sigma \right) &\xlongequal[]{\xi = 0}& \delta_l^0  \sqrt{\frac{1}{2\pi}} \frac{1}{\sigma^3}  \e^{-\sigma^{-2} \left (\rmvh^2 + \iota^2  \right)} ~.  \label{Kingxi0}
\end{eqnarray}
Here, $\delta_l^0$ is the Kronecker symbol.

\section{Normalized FPRS collision operator}
\label{Normalized FPRS collision operator}

Before presenting the expanded form of the normalized FPRS collision operator, we firstly calculate the derivatives of amplitudes of normalized Rosenbluth potential functions with respect to the speed coordinate $\vvbth$. The partial derivatives in axisymmetric velocity space can be formulated as:
  \begin{eqnarray}
      \ddscrvh \HhL \left(\vvbth,t \right) &=& \frac{1}{2 L + 1} \frac{1}{\vvbth^2} \left [-(L+1) \ILFhL + (L) \JLpFhL \right], \label{dvHhL}
      \\
      \ddscrvh \GhL \left(\vvbth,t \right) &=& \frac{(L-1) \ILFhL - (L) \JLnFhL}{(2 L -1) (2 L +1 )} - \frac{(L+1) \ILIIFhL - (L+2) \JLpFhL}{(2 L + 1) (2 L + 3 )} ~. \label{dvGhL}
  \end{eqnarray}
  Similarly, the second partial derivative of $\GhL$ with respect to $\vvbth$ is:
  \begin{eqnarray}
      \dddscrvh \GhL \left(\vvbth,t \right) &=& {C_{G_L}^{n}} \left(\ILFhL + \JLnFhL \right) + {C_{G_L}^{p}} \left(\ILIIFhL + \JLpFhL \right)~. \label{ddvGhL}
  \end{eqnarray}
  Here, the coefficients ${C_{G_L}^{n}}$ and ${C_{G_L}^{p}}$ are given by:
  \begin{eqnarray}
      {C_{G_L}^{n}} &= & - \frac{L(L-1)}{(2L-1) (2L+1)},  \label{CG2L} 
      \quad
      {C_{G_L}^{p}} \ = \ \frac{(L+1) (L+2)}{(2L+1) (2L+3)}  ~. \label{CG2L2}
  \end{eqnarray}
  Similarly, analytical expressions for the first two derivatives of $\fhl(\rmvh,t)$ with respect to $\rmvh$ can be derived based on Eq.\EQ{KFE}.

  Substituting the distribution function\EQ{fhfhl}, Eq.\EQ{FhFhL} and Rosenbluth potentials\EQ{HhHhL} into the normalized FPRS collision operator\EQ{FPShdab}, simplifying the result by combining the like terms yields the normalized FPRS collision operator in axisymmetric velocity space:
  \begin{eqnarray}
      \colhab \left(\vh,t \right) &=& 4 \pi\sum_{i=0}^{9} \Sh_i ~. \label{FPShdabSi}
  \end{eqnarray}
  The zero-order effect term resulting from  the normalized background distribution function, $\Fh$, in the collision term can be formulated as:
  \begin{eqnarray}
      \Sh_0 &=& m_M \sumLoLmax \FhL (\vvbth, t) \PL \times \sumlolmax \fhl (\rmvh, t) \Pl ~. \label{S0FAxis} 
  \end{eqnarray}
  The first-order effect terms resulting from $\HhL$ will be:
  \begin{eqnarray}
      \Sh_1 &=& \CHh \sumLoLmax \PL \ddscrvh \HhL \times \sumlolmax \Pl \ddrmvh \fhl , \label{S1HAxis} \\
      \Sh_2 &=& \CHh \frac{1}{\vvbth} \sumLILMax \HhL \PLI \times \frac{1}{\rmvh} \sumlIlMax \fhl \PlI ~. \label{S2HAxis} 
  \end{eqnarray}
  Another first-order effect term, $\Sh_3$, will be zero due to the azimuthal derivatives is zero. Specifically, when $m_a = m_b$, the coefficient $\CHh$ 
  provided in Eq.\EQ{CHh} will be zero, resulting in $\Sh_1=\Sh_2 \equiv 0$ as well. Similarly, the second-order effect terms related to $\GhL$ will be:
  \begin{eqnarray}
      \Sh_4 &=& \CGh \sumLoLmax \PL \dddscrvh \GhL \times \sumlolmax \Pl \dddrmvh \fhl , \label{S4GAxis} 
      \\
      \Sh_5 &=& 2 \CGh \frac{1}{\vvbth} \sumLILMax \left(\frac{\GhL}{\vvbth} - \PL \ddscrvh \GhL \right) 
      \times 
      \frac{1}{\rmvh}  \sumlIlMax \left(\frac{\fhl}{\rmvh} - \Pl \ddrmvh \fhl \right) \PlI, \label{S5GAxis} 
      \\
      \Sh_7 &=& \CGh \frac{1}{\vvbth} \sumLoLmax \left(\PLILII \frac{\GhL}{\vvbth} + \PL \dGhL \right) 
      \times 
      \frac{1}{\rmvh}  \sumlolmax \left(\PlIlII \frac{\fhl}{\rmvh} + \Pl \dfhl \right),  \label{S7GAxis} 
      \\
      \Sh_9 &=& \CGh \sumLoLmax \left(\PLIcotQ \frac{\GhL}{\vvbth^2} + \PL \frac{1}{\vvbth} \dGhL \right) 
      \times  
      \sumlolmax \left(\PlIcotQ  \frac{\fhl}{\rmvh^2} + \Pl \frac{1}{\rmvh} \dfhl \right),  \label{S9GAxis} 
  \end{eqnarray}
  where the coefficient $\CGh$ provided in Eq.\EQ{CGh} and coefficients
  \begin{eqnarray}
      \PlIcotQ &=&  \frac{\mu}{\sinQ}  \PlI, 
      \quad
      \PlIlII \ = \ \PlII  + \PlIcotQ ~.
  \end{eqnarray}
  The remaining second-order effect terms vanish, $\Sh_6=\Sh_8=0$. As demonstrated by equations above, the FPRS collision operator\EQ{FPShdabSi} is generally a nonlinear model.
  
  In particular, when the system is spherical symmetric in velocity space, the collision effect becomes independent of the angular coordinate of velocity space. Therefore, the normalized FPRS collision operator\EQgive{FPShdabSi} will be:
  \begin{eqnarray}
      \Sh_0 \left(\vh,t \right) &=& m_M \Fho (\vvbth, t) \times  \fho (\rmvh, t) , \label{S0F1} 
      \\
      \Sh_1 \left(\vh,t \right) &=& \CHh \ddscrvh \Hho \times \ddrmvh \fho , \label{S1H1} 
      \\ 
      \Sh_4 \left(\vh,t \right) &=& \CGh \dddscrvh \Gho \times \dddrmvh \fho , \label{S4G1}
      \\
      \Sh_7 \left(\vh,t \right) &=&  \Sh_9 = \CGh \frac{1}{\vvbth}  \ddscrvh \Gho \times \frac{1}{\rmvh}  \ddrmvh \fho ~. \label{S7G1} 
  \end{eqnarray}
  The remaining first-order effect terms and second-order effect terms will be zero, $\Sh_2=\Sh_3=\Sh_5=\Sh_6=\Sh_8=0$. 

  For self-collision process, the mass ratio $m_M \equiv 1$ and thermal velocity ratio $\vabth \equiv 1$. Therefore, the coefficients\EQgive{CGh} will remain constants, reads:
  \begin{eqnarray}
      \CHh = 0, \quad 
      \CGh = 1 / 2
      ~. \label{CFHGhaa}          
  \end{eqnarray}
  The normalized FPRS collision operator (give in Eqs.\EQ{FPShdabSi}-\EQo{S7G1}) can be further simplify to a self-collision operator by substituting $\FhL$ and $\vvbth$ with $\fhl$ and $\rmvh$ respectively.

\end{appendices}

\end{spacing}

\begin{spacing}{0.5}  

\begin{small}
\bibliography{Plasma}     

\bibliographystyle{ieeetr}   

\end{small}


\end{spacing}

\end{document}